\documentclass[12pt]{iopart}

\usepackage{iopams}
\usepackage{amssymb}
\usepackage{graphicx}
\usepackage{amsmath}
\newtheorem{theorem}{Theorem}

\newtheorem{corollary}[theorem]{Corollary}

\newtheorem{definition}[theorem]{Definition}

\newtheorem{lemma}[theorem]{Lemma}

\newtheorem{proposition}[theorem]{Proposition}
\newtheorem{remark}[theorem]{Remark}

\newenvironment{proof}[1][Proof]{\textbf{#1.} }{\ \rule{0.5em}{0.5em}}

\begin{document}

\title[Transition Tori in the Planar Restricted Elliptic Three Body Problem]{Transition Tori in the Planar Restricted Elliptic Three Body Problem}

\author{Maciej J. Capi\'nski}

\address{Faculty of Applied Mathematics, AGH University of Science and Technology,\\
al. Mickiewicza 30, 30-059 Krak\'ow, Poland.}
\ead{mcapinsk@agh.edu.pl}

\author{Piotr Zgliczy\'nski}

\address{Institute of Computer Science, Jagiellonian University, \\
Lojasiewicza 6, 30--348  Krak\'ow, Poland.}
\ead{piotr.zgliczynski@ii.uj.edu.pl}

\begin{abstract}
    We consider the elliptic three body problem as a perturbation of the circular problem.
    We show that for sufficiently small eccentricities of the elliptic problem, and
    for energies sufficiently close to the energy of the libration point $L_2$, a Cantor
    set of Lyapounov orbits survives the perturbation. The orbits are perturbed to quasi-periodic
    invariant tori. We show that for a certain family of masses of the primaries, for such tori
    we have transversal intersections of stable and unstable manifolds, which lead to chaotic
    dynamics involving diffusion over a short range of energy levels.
    Some parts of our argument are nonrigorous, but are strongly backed by numerical
    computations.
\end{abstract}

\footnotetext[1]{The research of the first author was partially
supported by the Polish Ministry of Science and Higher Education.
} \footnotetext[2]{The second author was supported in part by
Polish State Ministry of
            Science and Information Technology  grant N201 024 31/2163}
\footnotetext[1]{Both authors are supported by the Polish State Ministry of Science
and Information Technology grant N201 543238.}

\today

\maketitle

\section{Introduction}

In the planar restricted circular three body problem (PRC3BP) two large masses
$\mu$ and $1-\mu$ rotate on planar circular Keplerian orbits. For convenience
we will call the larger body - the Sun and the smaller massive body - the
planet The problem deals with the motion of a third massless particle (the
comet or the spacecraft), which moves on the same plane as the two larger
bodies under their gravitational pull. This problem was considered by Llibre,
Martinez and Simo in \cite{Simo} for energies of solutions close to the energy
of the libration point $L_{2}^{\mu}.$ There it has been shown that there
exists a family of parameters $\{\mu_{k}\}_{k=2}^{\infty}$ for which we have a
homoclinic orbit to the libration point $L_{2}^{\mu_{k}}$. Moreover, it has
been shown that for $\mu$ close to any of the values $\mu_{k},$ for a
Lyapounov orbit around $L_{2}^{\mu}$ with energy sufficiently close to the
energy of $L_{2}^{\mu}$, its stable and unstable manifolds intersect
transversally. This dynamics is restricted to a constant energy manifold and
leads to a homoclinic tangle and symbolic dynamics. Later a similar problem
has been numerically investigated by Koon, Lo, Marsden and Ross \cite{Marsden}%
, where smaller energies were considered. In such a case the chaotic dynamics
is extended to include homoclinic and heteroclinic tangles along stable and
unstable manifolds of Lyapounov orbits around both the libration point
$L_{1}^{\mu}$ and $L_{2}^{\mu}$. This has been later proven by Wilczak and
Zgliczy\'{n}ski using a method of covering relations and rigorous computer
assisted computations in \cite{Zgliczyn-Wilczak, Zgliczyn-Wilczak2}, for the
case of the Sun-Jupiter system and the energy of the comet Oterma.

All of the above mentioned results have a common feature: since the problem
follows from an autonomous Hamiltonian, the transversality of the
intersections and the chaotic dynamics of the system are always restricted to
a constant energy manifold. In this paper we are going to consider the planar
restricted elliptic three body problem (PRE3BP), where the equations are no
longer autonomous, which means that a change of energy of solutions is
possible. We will consider the circular problem considered by Llibre, Martinez
and Simo in \cite{Simo} and generalize it to allow the orbits of the planet
and the Sun to be elliptic with small eccentricities $e$. We will treat this
as a perturbation of the circular case. We will show that most of the
Lyapounov orbits around $L_{2}^{\mu}$ persist under such perturbation as
KAM-tori. Moreover, we will show that the symbolic dynamics associated with
these orbits also survives. This kind of the 'structural stability' of
symbolic dynamics constitute the main result of this paper. It will turn out
that we also have chaotic diffusion along the energy level. In effect the
dynamics of the elliptic problem is by one dimension richer than the dynamics
of the circular problem, where all solutions are restricted to a constant
energy manifold.

The diffusion between energies discussed in this paper follows from a
mechanism similar to the 1964 Arnold's example \cite{Arnold-arnold-diff}.
Arnold conjectured that this phenomenon appears in the three body problem. The
result of this paper is a small step towards a proof of this conjecture, but
the described dynamics does not fulfill all requirements. First of all, prior
to perturbation we do not have a fully integrable system. We start with the
circular problem with a setting in which we already have a transversal
homoclinic connection between Lyapounov orbits, which in our setting play the
role of lower dimensional normally hyperbolic invariant tori. Such systems are
referred to as "a priory unstable", or even "a priory chaotic". Secondly, and
most importantly, our diffusion is between energies with distance of order
$\left(  e\mu^{1/3}\right)  ^{1/2}$ and not of order one.

Throughout some of the so far explored examples a certain pattern can be
observed in the methods used for problems involving diffusion in priory
unstable systems (see for example \cite{Invariant-Tori} for a result of
Moeckel on detection of transition tori in the case of the planar five body
problem; or the work of Delshams, Llave and Seara \cite{Llave} on diffusion of
energy for perturbations of geodesic flows on a two dimensional torus; also
Wiggins \cite{Wiggins}, \cite{Wiggins-chaotic-transport} discusses this
mechanism in the case of perturbations of completely integrable systems).
First a normally hyperbolic invariant manifold foliated by invariant tori is
found. The tori are required to have hyperbolic stable and unstable manifolds
and a transversal intersection of these manifolds. Secondly a perturbation of
the system is considered. By perturbation theory (\cite{Herman},
\cite{Wiggins-norm-hyp}) of normally hyperbolic manifolds, the normally
hyperbolic invariant manifold and its stable and unstable manifolds persist
under the perturbation. Next step is to show that on the perturbed invariant
manifold most of the invariant tori survive. This under appropriate
nondegeneracy conditions is a result of the Kolmogorov Arnold Moser Theory
(KAM) \cite{Arnold-proof-kolmogorov},\cite{Kolmogorov}. Using KAM-technics
(for example \cite{Graff}, \cite{Herman} or \cite{Zehnder}) it can be shown
that most of the invariant tori persist and form a Cantor set having a
positive measure in the invariant manifold. The last step is to show that the
stable and unstable manifolds of the surviving tori intersect transversally.
This can be done by the use of a Melnikov type method along a homoclinic orbit
of the unperturbed problem. The transversal intersections between the
invariant manifolds of the perturbed tori lead to homoclinic tangles for each
of the surviving tori. In addition to this we also have a chaotic diffusion
along the Cantor set of homoclinic tangles between the tori. In this paper we
will follow this procedure.

When applying the method to prove the existence of transition chains for a
given physical problem the steps of the above described procedure, which
present the biggest obstacles are usually the verification of the assumptions
of the KAM theorem and computation of the Melnikov integral.

The fundamental role for our investigation is played by the Hill's
problem. In the neighborhood of the libration point the PCR3BP and
PER3BP, when written in a suitably rescaled Hill's coordinates,
are perturbations of the Hill's problem depending on two small
parameters $\mu$ and $e$. This gives us a 'local' picture around
$L_{2}^{\mu}$, the existence of normally hyperbolic invariant
manifold and KAM-tori. In our case the twist property needed for
the KAM theorem follows from the Lyapounov-Moser Theorem
\cite{Moser}. For sufficiently small $\mu$ we will prove the twist
property for the family of periodic orbits around $L_{2}^{\mu}$,
by approximating the PRC3BP with the Hill's problem and thus
obtain the following theorem (for a detailed formulation see
Theorem \ref{th:KAM for L2}).

\begin{theorem}
There exist positive constants $R_{\text{Hill}},\kappa,\mu^{\ast}%
\in\mathbb{R}$ such that for all mass parameters $\mu<\mu^{\ast}$
and any perturbation $e$ such that $e\mu^{-2/3}<\kappa$ most
Lyapounov orbits with radii in Hill's coordinates not exceeding
$R_{\text{Hill}}$ are perturbed to quasi periodic orbits (in other
words, to invariant two dimensional invariant tori in extended
phase space). The set of radiuses for which the tori survive forms
a Cantor set with complement measure smaller than $O(\left(  e\mu
^{-1/3}\right)  ^{1/2}).$
\end{theorem}

Our second result concerns the study of the stable and unstable
leaves of the Cantor set of surviving KAM tori. We do so by
applying a modification of a Melnikov method to obtain the
following result (for a detailed formulation see Theorem
\ref{th:Main-Detail}).

\begin{theorem}
\label{thm:main}Assume that for the sequence of masses $\{\mu
_{k}\}_{k=2}^{\infty}$ (the sequence is specified in Theorem
\ref{lem homoclinic orbit for uk}) the twist condition holds and
that a derivative of a Melnikov integral (\ref{eq:Melnikov}) at
zero is nonzero, then for any given $\mu_{k}$ there exists a
radius $R(\mu_{k})$ such that for perturbations $e$ with
$e\mu_{k}^{-2/3}<\kappa$ and sufficiently small $e\mu_{k}^{-1/3}$
there exists a homoclinic and a heteroclinic tangle between the
surviving tori or radii smaller than $R(\mu_{k})$. The tangle
implies existence of symbolic dynamics and diffusion in energy.
The diffusion occurs between surviving tori on an interval of
energies of order $\left(  e\mu _{k}^{-1/3}\right)  ^{1/2}.$
\end{theorem}

To apply Theorem \ref{thm:main} we need to back the argument by
numerical verification of its assumptions. We have verified the
twist condition numerically for large masses $\mu_{k}$ and
provided a rigorous argument that for sufficiently small masses
the condition holds true. For the Melnikov integral
(\ref{eq:Melnikov}) we can rigorously prove that it is convergent
and that it is zero at zero. It needs to be stressed though that
the assumption that the derivative of the Melnikov integral at
zero is nonzero has only been verified numerically.

We believe that the above mentioned numerical computations can be
performed using an rigorous-computer-assisted approach in the
spirit of \cite{Zgliczyn-Wilczak}. Such arguments require careful
estimates, use of topological tools, and are the subject of
ongoing work.

In our work we have been unable to obtain uniform bounds for the
size of the radii $R(\mu_{k})$ from Theorem \ref{thm:main} with
respect to $\mu _{k}$. From our proof it only follows that these
need to decrease together with $\mu_{k}$ so that at least
$\mu_{k}^{-1/3}R(\mu_{k})$ is smaller than some constant. It is
possible though that in a number of needed estimates $R(\mu_{k})$
has to be chosen even smaller. We remark also that the symbolic
dynamics proved in Theorem \ref{thm:main} will hold not only for
the family of parameters $\{ \mu_{k} \}$, but also for other
masses $\mu$ for which $| \mu_{k}-\mu| < \epsilon_{k} $ with
sufficiently small $\epsilon_{k}$. 

The paper is organized as follows. Section~\ref{sec:prelim} contains
preliminaries, where we recall the earlier results on the planar restricted
three body problem of \cite{Simo}, and introduce basic facts about the Hill's
problem and the PRE3BP. In Section~\ref{sec:lap-twist} we present the
Lyapounov--Moser theorem \cite{Moser} and show how to apply it to obtain the
twist property. In Section~\ref{sec:twist-L2} we show that we have a twist on
the family of Lyapounov orbits around $L_{2}^{\mu}$. In
Section~\ref{sec:norm-hyp-KAM} we apply the normally hyperbolic invariant
manifold theorem together with the KAM Theorem to show that most of the
Lyapounov orbits around $L_{2}^{\mu}$ persist under perturbation from the
circular problem to the elliptic problem and prove the first of our two main
theorems. In Section~\ref{sec:melnikov} we use a Melnikov type argument to
detect the transversal intersections between the stable and unstable manifolds
of the perturbed Lyapounov orbits. In Section~\ref{sec:meln-comp} we compute
the Melnikov integral. In Section~\ref{sec:tranchains} we gather together our
results and prove Theorem \ref{thm:main}.

\section{Preliminaries}

\label{sec:prelim}

\subsection{The Planar Restricted Circular Three Body Problem}

\label{subsec:PRC3BP}

In the planar restricted circular three body problem (PRC3BP) we consider the
motion of a small massless particle (a comet or a spacecraft), under the
gravitational pull of two larger bodies of mass $\mu$ and $1-\mu$ (called the
planet and the Sun, respectively) which move around the origin on circular
orbits of period $2\pi$ on the same plane as the massless body. The
Hamiltonian of the problem is given by \cite{Moser-ksiazka}
\begin{equation}
H(\mu,q,p,t)=\frac{p_{1}^{2}+p_{2}^{2}}{2}-\frac{1-\mu}{r_{1}(t)}-\frac{\mu
}{r_{2}(t)}, \label{eq: H for circular problem}%
\end{equation}
where $\left(  p,q\right)  =\left(  q_{1},q_{2},p_{1},p_{2}\right)  $ are the
coordinates of the massless particle and $r_{1}(t)$ and $r_{2}(t)$ are the
distances from the masses $1-\mu$ and $\mu$ respectively. After introducing a
new coordinates system $(x,y,p_{x},p_{y})$
\begin{equation}%
\begin{array}
[c]{ll}%
x=q_{1}\cos t+q_{2}\sin t, & \quad p_{x}=p_{1}\cos t+p_{2}\sin t,\\
y=-q_{1}\sin t+q_{2}\cos t, & \quad p_{y}=-p_{1}\sin t+p_{2}\cos t,
\end{array}
\label{eq:x,y-coordinates}%
\end{equation}
which rotates together with the two larger masses, the larger masses become
motionless and one obtains \cite{Moser-ksiazka} an autonomous Hamiltonian
\begin{equation}
H(\mu,x,y,p_{x},p_{y})=\frac{(p_{x}+y)^{2}+(p_{y}-x)^{2}}{2}-\Omega(x,y),
\label{eq:H-PRC3BP}%
\end{equation}
where
\begin{align}
\Omega(x,y)  &  =\frac{x^{2}+y^{2}}{2}+\frac{1-\mu}{r_{1}}+\frac{\mu}{r_{2}%
},\nonumber\\
r_{1}  &  =\sqrt{(x-\mu)^{2}+y^{2}},\quad r_{2}=\sqrt{(x+1-\mu)^{2}+y^{2}}.
\label{eq:r1-r2-prc3bp}%
\end{align}
The motion of the particle is given by the equation
\begin{equation}
\mathbf{\dot{x}}=J\nabla H(\mu,\mathbf{x}), \label{eq:PRC3BP}%
\end{equation}
where $\mathbf{x}=(x,y,p_{x},p_{y})\in\mathbb{R}^{4}$, $J=\left(
\begin{array}
[c]{cc}%
0 & Id\\
-Id & 0
\end{array}
\right)  $ and $Id$ is a two dimensional identity matrix.

\begin{figure}[h]
\begin{center}
\includegraphics[
height=1.1in,
width=4in
]{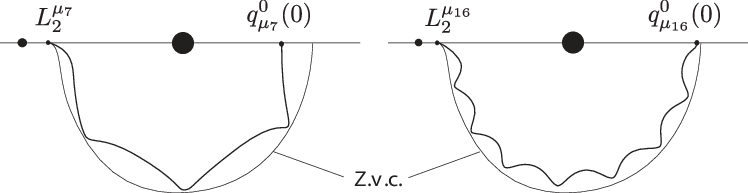}
\end{center}
\caption{The unstable manifold of $L_{2}^{\mu}$ in the $x,y$ coordinates.}%
\label{fig:wul2}%
\end{figure}

The movement of the flow (\ref{eq:PRC3BP}) is restricted to the hypersurfaces
determined by the energy level $C,$%
\begin{equation}
M(\mu,C)=\{(x,y,p_{x},p_{y})\in\mathbb{R}^{4}|H(\mu,x,y,p_{x},p_{y})=C\}.
\end{equation}
In the $x,y$ coordinates this means that the movement is restricted to the so
called Hill's region defined by
\begin{equation}
R(\mu,C)=\{(x,y)\in\mathbb{R}^{2}|\Omega(x,y)\geq-C\}.\nonumber
\end{equation}
\begin{figure}[h]
\begin{center}
\includegraphics[
height=1.4079in,
width=4.3907in
]{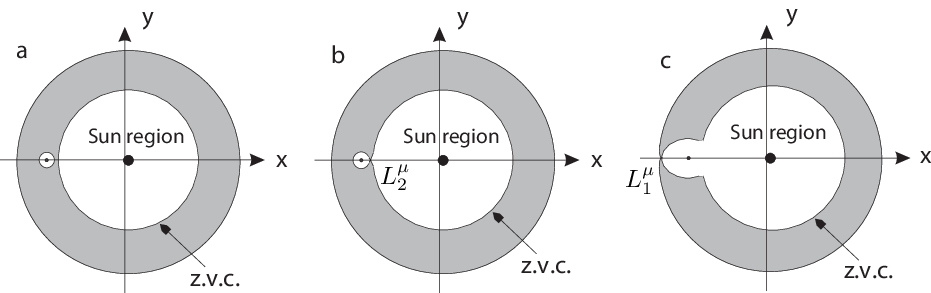}
\end{center}
\caption{The Hill's region for various energy levels: when the energy $C$ is
smaller than $C_{2}^{\mu}$ (a), when $C=C_{2}^{\mu}$ (b) and when
$C=C_{1}^{\mu}>C_{2}^{\mu}$ (c). }%
\label{fig:forbidden-region}%
\end{figure}The shape of the Hill's region $R(\mu,C)$ will differ with $C$
(see Figure \ref{fig:forbidden-region}). The focus of our attention in this
paper will be on the case when the energy $C$ is equal to or slightly larger
than $C_{2}^{\mu}.$ For the energy $C$ equal to $C_{2}^{\mu}$ we have the
libration point $L_{2}^{\mu}$ which is of the form $(-k,0,0,-k)$ with $k>0$.
We shall investigate the dynamics inside of the inner part of the Hill's
region to the right of the point $L_{2}^{\mu}$. We shall refer to it as the
"Sun region" (see Figure \ref{fig:forbidden-region}). The boundary of this
region (see Figures \ref{fig:forbidden-region} and \ref{fig:wul2}) is a zero
velocity curve (z.v.c.). The linearized vector field at the point $L_{2}^{\mu
}$ has two real and two purely imaginary eigenvalues, thus it follows
\cite{Simo} from the Lyapounov theorem that for energies $C$ larger and
sufficiently close to $C_{2}^{\mu}$ there exists a family of periodic orbits
$l_{\mu}(C)$ emanating from the equilibrium point $L_{2}^{\mu}.$

The PRC3BP admits the following reversing symmetry
\begin{equation}
S(x,y,p_{x},p_{y})=(x,-y,-p_{x},p_{y}). \label{eq:S-def}%
\end{equation}
We will say that an orbit $q(t)$ is $S$-symmetric when
\begin{equation}
S(q(t))=q(-t). \label{eq:S-sym}%
\end{equation}
In PRC3BP the Lyapounov orbits are $S$-symmetric (we have to choose the
initial time so that orbits start from the section $\{y=0\}$ at time $t=0)$.

We have the following results about the stable and unstable manifolds of
$L_{2}^{\mu}$ and $l_{\mu}(C)$.

\begin{theorem}
[{\cite[Theorem A]{Simo}}]\label{lem homoclinic orbit for uk} For $\mu$
sufficiently small the branch of $W_{L_{2}^{\mu}}^{u}$ contained in the Sun
region (see Figures \ref{fig:forbidden-region} and \ref{fig:wul2}) has a
projection on the bounded component of $R(\mu,C)$ given by
\begin{align}
d(t)  &  =\mu^{1/3}\left(  \frac{2}{3}N(\infty)-3^{1/6}+M(\infty)\cos
t+o(1)\right)  ,\label{eq: d=}\\
\alpha(t)  &  =-\pi+\mu^{1/3}\left(  N(\infty)t+2M(\infty)\sin t+o(1)\right)
, \label{eq: alpha=}%
\end{align}
where $d$ is the distance to the z.v.c., $\alpha$ the angular coordinate,
$N(\infty)$ and $M(\infty)$ are constants and the expressions remain true out
of a given neighborhood of $L_{2}^{\mu}.$ The parameter $t$ means the physical
time from a suitable origin. The terms $o(1)$ tend to zero when $\mu$ does and
they are uniform in $t$ for $t=O(\mu^{-1/3}).$

In particular the first intersection with the $x$ axis is orthogonal to that
axis, giving a $S$-symmetric homoclinic orbit for a sequence of values $\mu$
which has the following asymptotic expression:
\begin{equation}
\mu_{k}=\frac{1}{N(\infty)^{3}k^{3}}(1+o(1)). \label{eq:muk}%
\end{equation}

\end{theorem}

Let us now introduce a notation for the $S$-symmetric homoclinic orbit to
$L_{2}^{\mu_{k}}$ obtained in Theorem \ref{lem homoclinic orbit for uk} for
the parameters $\mu_{k}$ given in (\ref{eq:muk}). We will denote such an orbit
by $q_{\mu_{k}}^{0}(t)$ (see Figure \ref{fig:wul2}). We assume that such an
orbit starts at a section $\{y=0\}$ at time $t=0$.

\begin{figure}[h]
\begin{center}
\includegraphics[
height=2.1395in,
width=2.4483in
]{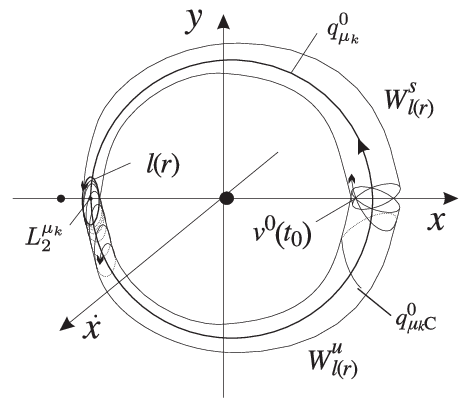}
\end{center}
\caption{The intersections of the stable and unstable manifolds of $l_{\mu
}\left(  C\right)  $ in the PRC3BP.}%
\label{fig:tubes}%
\end{figure}

\begin{theorem}
[{\cite[Theorem B]{Simo}}]\label{lem: transversality-Simo}For $\mu$ and
$\Delta C=C-C_{2}^{\mu}$ sufficiently small, the branch of $W^{u}\left(
l_{\mu}\left(  C\right)  \right)  $ contained in the Sun region intersects the
plane $y=0$ for $x>0$ in a curve diffeomorphic to a circle (see Figure
\ref{fig:tubes}) given by%
\begin{align*}
x  &  =x_{w}-\sqrt{\Delta C}\left(  N+2M\cos\tau\right)  ^{-1}\left(  2M+N\cos
M_{f}\right)  (K_{1}\cos\tau\cos\sigma-K_{2}\sin\tau\sin\sigma)\\
&  +\mu^{1/3}M(1-\cos M_{f})\\
&  +\mu^{2/3}\left\{  -\frac{2MN}{3}\left(  1-\cos M_{f}\right)  +M^{2}%
\sin^{2}M_{f}-\frac{2}{9}N\alpha-\frac{M}{3}\alpha\cos M_{f}\right\}
+O(\mu),\\
\dot{x}  &  =\sin M_{f}[\sqrt{\Delta C}N\left(  N+2M\cos\tau\right)
^{-1}\left(  K_{1}\cos\tau\cos\sigma-K_{2}\sin\tau\sin\sigma\right) \\
&  +\mu^{1/3}M+\mu^{2/3}\left\{  \frac{MN}{3}+2M^{2}\cos M_{f}+\frac{M}%
{3}\alpha\right\}  ]+O(\mu)
\end{align*}
where $x_{w},$ $M,$ $N,$ $\tau,$ $K_{1},$ $K_{2}$ are suitable constants,
$\alpha$ measures the distance from $\mu$ to some $\mu_{k}$ given by Theorem
\ref{lem homoclinic orbit for uk}, $M_{f}$ is obtained implicitly from%
\begin{gather}
\left(  \frac{1}{3}N\alpha+\mu^{-2/3}\sqrt{\Delta C}3M\left(  N+2M\cos
\tau\right)  ^{-1}\left(  K_{1}\cos\tau\cos\sigma-K_{2}\sin\tau\sin
\sigma\right)  \right)  \frac{\pi}{N}\nonumber\\
+NM_{f}+2M\sin M_{f}=o(1), \label{eq:Mf-implicit}%
\end{gather}
and $\sigma,$ ranging from $0$ to $2\pi,$ is the parameter of the curve.

Moreover for points in the $(\mu,C)$ plane such that there exists a $\mu_{k}$
of Theorem \ref{lem homoclinic orbit for uk} for which
\begin{equation}
\Delta C>L\mu_{k}^{4/3}(\mu-\mu_{k})^{2} \label{eq:DeltaC-L}%
\end{equation}
holds (where $L$ is a constant), there exist $S$-symmetric transversal
homoclinic orbits. In particular, for $\mu=\mu_{k}$ there exist symmetrical
transversal homoclinic orbits $q_{\mu_{k}C}^{0}$ for the periodic orbit
$l_{\mu_{k}}(C)$ for $C>C_{2}^{\mu_{k}}$ arbitrarily close to $C_{2}^{\mu_{k}%
}.$
\end{theorem}

\begin{remark}
In the original version of Theorems \ref{lem homoclinic orbit for uk} and
\ref{lem: transversality-Simo} the $C$ was taken as the Jacobi constant $C=F$
where
\begin{equation}
F(x,y,p_{x},p_{y})=-2H(\mu,x,y,p_{x},p_{y})=2\Omega(x,y)-\left(  \dot{x}%
^{2}+\dot{y}^{2}\right)  . \label{eq:Jacobi-int}%
\end{equation}
In this paper we have rewritten the Theorems with $C$ as the Hamiltonian of
the PRC3BP, which means that we have a change of sign in $C$ compared with the
original version.
\end{remark}

\begin{remark}
\label{rem:symb-dyn} Using a standard dynamical system theory argument, from
the Birkhoff-Smale homoclinic theorem, the transversal homoclinic connections
to Lyapounov orbits imply chaotic symbolic dynamics of the system. This is a
content of Theorem C in \cite{Simo}. Since the system is autonomous this
dynamics is restricted to the constant energy manifold.
\end{remark}

\begin{remark}
\label{rem:curve-splitting} For sufficiently small $\mu=\mu_{k}$, the curves
(obtained in Theorem \ref{lem: transversality-Simo}) on $\{y=0\}$ associated
with the stable an unstable manifolds of $l_{\mu_{k}}\left(  C\right)  $
intersect transversally at an angle $O(\mu_{k}^{1/3})$. This can be derived
based on the parameterization of the curves from Theorem
\ref{lem: transversality-Simo} combined with the symmetry property
(\ref{eq:S-sym}) of the PRC3BP. This is done in the Appendix. For more details
on the interpretation of the curves from the theorem see also \cite{Simo}.
\end{remark}

\subsection{The Hill's Problem}

\label{subsec:Hill}

Let us consider a change of coordinates which shifts the origin to the smaller
body of the mass $\mu$ and rescales the coordinates by the factor $\mu
^{-1/3}.$ We will refer to the following as the Hill's coordinates.%
\begin{equation}
\mathbf{\bar{x}}=\mu^{-1/3}\left(  \mathbf{x}-\left(  \mu-1,0,0,\mu-1\right)
\right)  . \label{eq:Hill-coordinates}%
\end{equation}
We will rewrite the Hamiltonian (\ref{eq:H-PRC3BP}) and derive the formula of
the Hill's problem and list a number of facts which will be relevant for us in
the future.

Let us start with a simple lemma.

\begin{lemma}
\label{lem:Hscaling} Consider a Hamiltonian $H:\mathbb{R}^{n}\times
\mathbb{R}^{n}\rightarrow\mathbb{R}$ and a transformation $\bar{p}=\beta p$,
$\bar{q}=\beta q$. Then $(p(t),q(t))$ is a solution of the Hamiltonian system
for $H(p,q)$ if and only if $(\bar{p}(t),\bar{q}(t))$ is a solution of the
Hamiltonian system for $\bar{H}(\bar{p},\bar{q})=\beta^{2}H\left(  \frac
{\bar{p}}{\beta},\frac{\bar{q}}{\beta}\right)  $.
\end{lemma}

The shift of the origin present in the transformation
(\ref{eq:Hill-coordinates}) is clearly canonical hence we can use Lemma
\ref{lem:Hscaling} to obtain the Hamiltonian in new variables with $\beta
=\mu^{-1/3}$%
\begin{equation}
\bar{H}(\mu,\mathbf{\bar{x}})=\mu^{-2/3}H\left(  \mu,\mu^{1/3}\mathbf{\bar{x}%
}+\left(  \mu-1,0,0,\mu-1\right)  \right)  . \label{eq:PRC3BP-Hill}%
\end{equation}
By expanding $\Omega$ in the new coordinates around zero we can rewrite our
family of Hamiltonians as
\begin{align}
\bar{H}(\mu,\mathbf{\bar{x}})  &  =\bar{H}(\mu,\bar{x},\bar{y},\bar{p_{x}%
},\bar{p_{y}})=\frac{(\bar{p}_{x}+\bar{y})^{2}+(\bar{p}_{y}-\bar{x})^{2}}%
{2}+\label{eq:PRC3BP-Hill-reduced}\\
&  -\frac{1}{\bar{r}}-\frac{3}{2}\bar{x}^{2}+O(\mu^{1/3})+C(\mu),\quad
\mbox{for $\bar{r} <\alpha
\mu^{-1/3}$}\nonumber
\end{align}
where $C(\mu)=\mu^{-2/3}\left(  \left(  1-\mu\right)  +\left(  1-\mu\right)
^{2}/2\right)  $ and $\bar{r}=\sqrt{\bar{x}^{2}+\bar{y}^{2}}$. The term
$O(\mu^{1/3})$ depends on $(\bar{x},\bar{y})$ and can be written as a function
$a(\mu,\bar{x},\bar{y})$, which for any sufficiently small $\alpha<1$ and
$\mu\in\lbrack0,1]$, $\bar{r}<\alpha\mu^{-1/3}$ satisfies
\begin{equation}
\frac{|a(\mu,\bar{x},\bar{y})|}{r^{3}}\leq M(\alpha).\nonumber
\end{equation}
The reason for introducing $\alpha$ is, that we have to be away from the Sun
in order for Taylor series of $\frac{1}{r_{1}}$ to be convergent. Observe that
we can drop the term $C(\mu)$.

It should be stressed that $\bar{H}$ depends analytically on $\mathbf{\bar{x}%
}$ and $\mu^{1/3}$, hence the derivatives of the $O(\mu^{1/3})$ with respect
to $\mathbf{\bar{x}}$ are still $O(\mu^{1/3})$. Therefore for fixed $\mu$
\begin{equation}
\mathbf{\bar{x}}^{\prime}=J\nabla\bar{H}(\mu,\mathbf{\bar{x}})=J\nabla\bar
{H}(0,\mathbf{\bar{x}})+O(\mu^{1/3}).\nonumber
\end{equation}
The term $O(\mu^{1/3})$ is uniform in $\mathbf{\bar{x}}$ for $|\mathbf{\bar
{x}}|\leq\alpha\mu^{-1/3}$ for any fixed sufficiently small $\alpha<1$. The
Hamiltonian $\bar{H}(0,\mathbf{\bar{x}})$ is the Hamiltonian of the Hill's
problem%
\begin{equation}
H^{\text{Hill}}(\mathbf{\bar{x}})=\bar{H}(0,\mathbf{\bar{x}}).
\label{eq: Hills problem}%
\end{equation}

If we denote by $q^{\text{Hill}}(t)$ the solution of the Hill problem and by
$q^{\mu}(t)$ the solution of the PCR3BP, both expressed in Hill's coordinates
(\ref{eq:Hill-coordinates}) and both starting from the same initial condition,
then the following holds
\begin{equation}
|q^{\text{Hill}}(t)-q^{\mu}(t)|\leq e^{l|t|}O(\mu^{1/3}),
\label{eq:estm-Hill-PCR3BP}%
\end{equation}
provided there exist $0<\alpha<1$ and a compact convex set $Z\subset
B(0,\alpha\mu^{-1/3})\times\mathbb{R}^{2}$, such that $Z$ contains both
$q^{\text{Hill}}([0,t])$ and $q^{\mu}([0,t])$ and $(0,0)\notin\pi_{x,y}(Z)$.
The constant $l$ depends on $Z$.

Let us now list a few properties of the Hill's problem. The problem has two
equilibrium points, $L_{1}^{\text{Hill}}=(-3^{-1/3},0,0,-3^{-1/3})$ and
$L_{2}^{\text{Hill}}=(3^{-1/3},0,0,3^{-1/3})$. The linearization of
$x^{\prime}=J\nabla H^{Hill}$ at $L_{2}^{Hill}$ is given by $x^{\prime}=Ax$,
where\emph{ }
\begin{equation}
A=\left(
\begin{array}
[c]{cccc}%
0 & 1 & 1 & 0\\
-1 & 0 & 0 & 1\\
8 & 0 & 0 & 1\\
0 & -4 & -1 & 0
\end{array}
\right)  . \label{eq:Hill-linear}%
\end{equation}
The eigenvalues of $A$ are: $\pm\alpha_{1},$ $\pm\alpha_{2}$ with $\alpha
_{1}=\sqrt{1+2\sqrt{7}}$ and $\alpha_{2}=\sqrt{1-2\sqrt{7}}.$ The Hill problem
has a reversing symmetry $S$ given by (\ref{eq:S-def}).

\subsection{The Planar Restricted Elliptic Three Body Problem}

The planar restricted elliptic three body problem (PRE3BP) differs from the
PRC3BP by the fact that the two larger bodies move on elliptic orbits of
eccentricities $e$ instead of circular orbits. The period of these orbits is
$2\pi$ and the Hamiltonian of the PRE3BP is analogous to
(\ref{eq: H for circular problem}), with the only difference that in
$r_{1}(t)$ and $r_{2}(t)$ we take the distance from the elliptic instead of
the circular orbits of the two larger masses. The trajectories of these orbits
can be written as (see \cite{Xia}) $((\mu-1)x(t),(\mu-1)y(t))$ for the body
$\mu$ and $\left(  \mu x\left(  t\right)  ,\mu y\left(  t\right)  \right)  $
for $\left(  1-\mu\right)  $, where%
\begin{align}
x(t)  &  =(1-e\cos\psi)\cos\psi+O(e^{2}),\nonumber\\
y(t)  &  =(1-e\cos\psi)\sin\psi+O(e^{2}),\label{eq:eliptic-orbits}\\
\psi(t)  &  =t+2e\sin t+O(e^{2}).\nonumber
\end{align}

If one changes into the rotating coordinates (\ref{eq:x,y-coordinates}), which
is a canonical transformation (see \cite{Meyer-Hall}), then the Hamiltonian
(\ref{eq: H for circular problem}) becomes (for a detailed derivation see the
Appendix)
\begin{equation}
H^{e}(\mu,\mathbf{x},t)=H(\mu,\mathbf{x})+eG(\mu,\mathbf{x},t)+O(e^{2}%
\mu^{-2/3}), \label{eq:H-PRE3BP}%
\end{equation}
where $H$ is the Hamiltonian of the PRC3BP (\ref{eq:H-PRC3BP}), $r_{2}$ is
given in (\ref{eq:r1-r2-prc3bp}), $G$ is $2\pi$ periodic over $t$ and is given
by the formula
\begin{equation}
G=\frac{1-\mu}{\left(  r_{1}\right)  ^{3}}\bar{g}(\mu,x,y,t)+\frac{\mu
}{\left(  r_{2}\right)  ^{3}}\bar{g}(\mu-1,x,y,t), \label{eq:Gdef}%
\end{equation}%
\begin{equation}
\bar{g}(\alpha,x,y,t)=\alpha(-2y\sin{t}+x\cos{t})-\alpha^{2}\cos{t}.
\label{eq:g-bar}%
\end{equation}
The term $O( e^{2}\mu^{-2/3})=a(x,y,e,\mu,t)$ satisfies $\frac{|a(x,y,e,\mu
,t)|}{e^{2}\mu^{-2/3}}<M(\delta,\kappa,R)$, on the set defined by the
following conditions
\begin{align}
\mu &  \in[0,1],\quad t\in\mathbb{R},\quad e\mu^{-2/3}<\kappa
\label{eq:kappa-e-bound}\\
r_{1}  &  >\delta,\quad r_{2}\geq\delta\mu^{1/3},\quad\sqrt{x^{2}+y^{2}}\leq
R,\nonumber
\end{align}
were $\delta>0$ measures the closest approach to the Sun and to the planet
(multiplied by $\mu^{-1/3}$) is a number around $1/2$, $R$ is the radius of
the ball containing all orbits of interest in our problem, hence we can take
$R=2$ and $\kappa$ is sufficiently small number (for more details on $\kappa$
see the derivation in the Appendix).

We consider points $(x,y)$ inside of the "Sun region" (see Figure
\ref{fig:forbidden-region}) which means that $r_{2}>\frac{1}{2}\left\Vert
L_{2}^{\mu}-\left(  \mu-1,0,0,\mu-1\right)  \right\Vert >\delta\mu^{1/3}$ for
some $\delta>0.$ The Hamiltonian (\ref{eq:H-PRE3BP}) can be rewritten in
Hill's coordinates (\ref{eq:Hill-coordinates}) as
\begin{equation}
\bar{H}^{e}(\mu,\mathbf{\bar{x}},t)=\bar{H}(\mu,\mathbf{\bar{x}})+e\bar
{G}\left(  \mu,\mathbf{\bar{x}},t\right)  +O(e^{2}\bar{r}^{2}),
\label{eq:PRE3BP-Hill}%
\end{equation}%
\begin{equation}
\bar{G}\left(  \mu,\mathbf{\bar{x}},t\right)  =\mu^{-2/3}G(\mu,\mu
^{1/3}\mathbf{\bar{x}+}\left(  \mu-1,0,0,\mu-1\right)  ,t), \label{eq:G-bar}%
\end{equation}

Further, in order to understand the mutual relation between $e$ and $\mu$ for
which we have interesting dynamical phenomena, it will be important to observe
that after rearranging and neglecting terms independent of $\bar{x}$,
(\ref{eq:PRE3BP-Hill}) can be written as
\begin{equation}
\bar{H}^{e}(\mu,\mathbf{\bar{x}},t)=\bar{H}(\mu,\mathbf{\bar{x}})+e\mu
^{-1/3}\frac{2\bar{y}\sin{t}-\bar{x}\cos{t}}{\bar{r}^{3}}+O(e\mu
^{2/3})+O(e^{2}\mu^{-4/3}), \label{eq:He-reduced}%
\end{equation}
on the following set
\begin{align*}
\mu &  \in\lbrack0,1],\quad t\in\mathbb{R},\quad e\mu^{-2/3}<\kappa\\
\bar{r}  &  >\delta,\quad\bar{r}\leq M_{1},\quad M_{1}\mu^{1/3}<1.
\end{align*}

The Hamiltonian $\bar{H}^{e}$ generates a differential equation
\begin{equation}
\mathbf{\bar{x}}^{\prime}=f(\mu,\mathbf{\bar{x}})+eg(\mu,\mathbf{\bar{x}%
},t)+O(e\mu^{1/3})+O(e^{2}\mu^{-4/3}) \label{eq:PRE3BP-reduced}%
\end{equation}
where
\begin{align}
f(\mu,\mathbf{x})  &  =J\nabla\bar{H}(\mu,\mathbf{x})\label{eq:def-of-f}\\
g(\mu,\mathbf{x},t)  &  =J\nabla\bar{G}(\mu,\mathbf{x}). \label{eq:def-of-g}%
\end{align}

\begin{remark}
\label{rem:range-rh} In our future consideration we shall use equation
(\ref{eq:He-reduced}) in a neighborhood of $\bar{L}_{2}^{\mu}=\mu
^{-1/3}\left(  L_{2}^{\mu}-\left(  \mu-1,0,0,\mu-1\right)  \right)  $ of
constant radius (later denoted as $R_{Hill}$), in which Lyapounov orbits
reside. In such neighborhood the term $O(e\mu^{2/3})$ of (\ref{eq:He-reduced})
together with higher order terms are uniform. It needs to be emphasized though
that in order for (\ref{eq:He-reduced}) to be valid for a given $\mu$, we
first need to choose $e$ sufficiently small so that the estimate $e\mu
^{-2/3}<\kappa$ in (\ref{eq:kappa-e-bound}) holds true.
\end{remark}

\section{From Lyapounov-Moser Theorem to Twist Property at Equilibrium Points}

\label{sec:lap-twist}

In this section we shall show how one can prove the twist property at an
equilibrium point using the Lyapounov-Moser Theorem \cite{Moser}. First the
Theorem will be stated. Next a number of observations on the Theorem in the
special case of one real and one pure imaginary eigenvalue with just two
degrees of freedom will be made. This will be followed by a brief outline of
the construction by which the Theorem was proved \cite{Moser} from which the
twist property will follow.

\begin{theorem}
[The Lyapounov-Moser Theorem \cite{Moser}]\label{th: Twierdzenie Mosera} Let
\begin{align}
\dot{q}_{\nu}  &  =H_{p_{\nu}}(q,p)\label{eq: x'=JGradH}\\
\dot{p}_{\nu}  &  =-H_{q_{\nu}}(q,p)\nonumber
\end{align}
$\nu=1,\ldots,n,$ be an analytic Hamiltonian system with $n$ degrees of
freedom and an equilibrium solution $q=p=0$. Let $\alpha_{1},\ldots,\alpha
_{n},-\alpha_{1},\ldots,-\alpha_{n}$ be the eigenvalues of the linearization
of (\ref{eq: x'=JGradH}) at the equilibrium point. Assume that the
eigenvalues
\begin{equation}
\alpha_{1},\ldots,\alpha_{n},-\alpha_{1},\ldots,-\alpha_{n}\nonumber
\end{equation}
are $2n$ different complex numbers and that $\alpha_{1},\alpha_{2}$ are
independent over the reals. Let us also assume that for any integer numbers
$n_{1}$ and $n_{2}$
\begin{equation}
\alpha_{\nu}\neq n_{1}\alpha_{1}+n_{2}\alpha_{2}\quad\text{for }\nu
\geq3.\nonumber
\end{equation}
Then there exists a four parameter family of solutions of (\ref{eq: x'=JGradH}%
) of the form
\begin{align}
q_{\nu}  &  =\phi_{\nu}(\xi_{1},\xi_{2},\eta_{1},\eta_{2}) \label{eq: x_v=phi}%
\\
p_{\nu}  &  =\psi_{\nu}(\xi_{1},\xi_{2},\eta_{1},\eta_{2})\nonumber
\end{align}
where
\begin{equation}
\xi_{k}(t)=\xi_{k}^{0}e^{ta_{k}(\xi_{1}^{0}\eta_{1}^{0},\xi_{2}^{0}\eta
_{2}^{0})},\quad\eta_{k}(t)=\eta_{k}^{0}e^{-ta_{k}(\xi_{1}^{0}\eta_{1}^{0}%
,\xi_{2}^{0}\eta_{2}^{0})}\quad\text{for }k=1,2, \label{eq: xi_k=...}%
\end{equation}
and
\begin{equation}
a_{1}(\xi_{1}^{0}\eta_{1}^{0},\xi_{2}^{0}\eta_{2}^{0})=\alpha_{1}+...,\qquad
a_{2}(\xi_{1}^{0}\eta_{1}^{0},\xi_{2}^{0}\eta_{2}^{0})=\alpha_{2}+...
\label{eq:a-series}%
\end{equation}
are convergent power series. The series $\phi_{\nu},\psi_{\nu}$ converge in
the neighborhood of the origin and the rank of the matrix
\begin{equation}
\left(
\begin{array}
[c]{cc}%
\phi_{\nu\xi_{k}} & \phi_{\nu\eta_{k}}\\
\psi_{\nu\xi_{k}} & \psi_{\nu\eta_{k}}%
\end{array}
\right)  _{\substack{\nu=1,2,\ldots,n\\k=1,2}}\nonumber
\end{equation}
is four. The solutions (\ref{eq: x_v=phi}) depend on four small enough complex
parameters $\xi_{k}^{0},\eta_{k}^{0}.$

If in addition $\alpha_{1}$, $\alpha_{2}$, $-\alpha_{1}$, $-\alpha_{2}$
contain their complex conjugates, the solution can be chosen to be real,
depending on 4 real parameters.
\end{theorem}

In the case of the PRC3BP, $n$ is simply equal to two and the equations
(\ref{eq: x_v=phi}), (\ref{eq: xi_k=...}) describe all the solutions near the
neighborhood of the equilibrium point. We will be interested in the
application of the Theorem to the libration point $L_{2}^{\mu},$ where
$\alpha_{1}$ is real and $\alpha_{2}$ is pure imaginary. From now on we will
restrict our discussion to this particular case. The following remarks and
lemmas adapt Theorem \ref{th: Twierdzenie Mosera} to this setting.

\begin{remark}
\label{lem: a1 real a2 imaginary} When the system (\ref{eq: x'=JGradH}) is
generated by a real Hamiltonian and if $\alpha_{1}$ is real and $\alpha_{2} $
is pure imaginary then for the real solutions of (\ref{eq: x'=JGradH}) of the
form (\ref{eq: x_v=phi}) the functions $\xi_{k}(t)$ and $\eta_{k}(t)$ are
invariant under the involution \cite[page 102]{Siegel}
\begin{equation}
J_{w}(\xi_{1},\xi_{2},\eta_{1},\eta_{2})=(\bar{\xi}_{1},i\bar{\eta}_{2}%
,\bar{\eta}_{1},i\bar{\xi}_{2}). \label{eq:def-inv}%
\end{equation}
Let us also note that the original version of Theorem
\ref{th: Twierdzenie Mosera} in \cite{Moser} contained an error. There an
involution
\begin{equation}
J_{w}(\xi_{1},\xi_{2},\eta_{1},\eta_{2})=(\bar{\xi}_{1},\bar{\eta}_{2}%
,\bar{\eta}_{1},\bar{\xi}_{2})\nonumber
\end{equation}
was proposed. This stands in conflict with a requirement that the
transformation $\Phi$ in the proof of Theorem \ref{th: Twierdzenie Mosera}
\cite{Moser} should be canonical (See also equations (\ref{eq:Phi-kan}) and
(\ref{eq: reality condition})).
\end{remark}

The reality condition (the fact that a point $(\xi_{1},\xi_{2},\eta_{1}%
,\eta_{2})$ given in new coordinates represents a point from $\mathbb{R}^{4}$
in original ones) is
\begin{equation}
J_{w}(\xi_{1},\xi_{2},\eta_{1},\eta_{2})=(\xi_{1},\xi_{2},\eta_{1},\eta_{2}).
\label{eq:reality}%
\end{equation}
The subspace of $\mathbb{C}^{4}$ of fixed points of $J_{\omega}$,
$\mbox{Fix}(J_{\omega})$ is given by
\begin{equation}
\xi_{1}\in\mathbb{R},\eta_{1}\in\mathbb{R},\xi_{2}=re^{i\varphi},\eta
_{2}=ire^{-i\varphi},\nonumber
\end{equation}
where $r,\varphi\in\mathbb{R}$. On $\mbox{Fix}(J_{w})$ we will use the
coordinates $(\xi_{1},\eta_{1},r,\varphi)$.

\begin{remark}
\label{rem-radius-convergence}From the proof of the convergence of the series
(\ref{eq: x_v=phi}), (\ref{eq:a-series}) during the proof of Theorem
\ref{th: Twierdzenie Mosera} in \cite{Moser}, it follows that if we consider a
family of Hamiltonians
\begin{equation}
H_{\lambda}:\mathbb{R}^{n}\times\mathbb{R}^{n}\rightarrow\mathbb{R},\nonumber
\end{equation}
which is an analytic function of all variables including the parameter and
possesses for each $\lambda\in I$, where $I$ is a closed interval, a (locally)
unique fixed point $p_{\lambda}$ of center-saddle type depending analytically
on $\lambda$, then the radius of convergence of the series (\ref{eq: x_v=phi}%
), (\ref{eq:a-series}) around $p_{\lambda}$ can be chosen uniformly for close
values of $\lambda$.
\end{remark}

\begin{lemma}
\label{lem: a2 imaginary} If $\alpha_{1}$ is real and $\alpha_{2}$ is pure
imaginary then for all real solutions of (\ref{eq: x'=JGradH}) the series
$a_{1}$ from the Theorem \ref{th: Twierdzenie Mosera} is real and the series
$a_{2}$ is pure imaginary. Moreover if we choose a real periodic solution
\begin{equation}%
\begin{array}
[c]{c}%
q_{\nu}(t)=\phi_{\nu}(0,\xi_{2}(t),0,\eta_{2}(t))\\
p_{\nu}(t)=\psi_{\nu}(0,\xi_{2}(t),0,\eta_{2}(t))
\end{array}
\quad\nu=1,2 \label{eq: xv, yv = ...}%
\end{equation}
where $\xi_{2}(t)$ and $\eta_{2}(t)$ are given by (\ref{eq: xi_k=...}), then
there exist two real numbers $r\ $and $\varphi$ such that
\begin{align*}
\xi_{2}(t)  &  =re^{ta_{2}(0,ir^{2})+i\varphi}\\
\eta_{2}(t)  &  =ire^{-ta_{2}(0,ir^{2})-i\varphi}.
\end{align*}

\end{lemma}

\begin{proof}
From Remark \ref{lem: a1 real a2 imaginary} we know that the real solutions
satisfy the reality condition (\ref{eq:reality}). We therefore have
\begin{align}
\xi_{1}^{0}e^{ta_{1}(\xi_{1}^{0}\eta_{1}^{0},\xi_{2}^{0}\eta_{2}^{0})}  &
=\overline{\xi_{1}^{0}e^{ta_{1}(\xi_{1}^{0}\eta_{1}^{0},\xi_{2}^{0}\eta
_{2}^{0})}}\label{eq:inv1}\\
\xi_{2}^{0}e^{ta_{2}(\xi_{1}^{0}\eta_{1}^{0},\xi_{2}^{0}\eta_{2}^{0})}  &
=i\left(  \overline{\eta_{2}^{0}e^{-ta_{2}(\xi_{1}^{0}\eta_{1}^{0},\xi_{2}%
^{0}\eta_{2}^{0})}}\right) \label{eq:inv2}\\
\eta_{1}^{0}e^{-ta_{1}(\xi_{1}^{0}\eta_{1}^{0},\xi_{2}^{0}\eta_{2}^{0})}  &
=\overline{\eta_{1}^{0}e^{-ta_{1}(\xi_{1}^{0}\eta_{1}^{0},\xi_{2}^{0}\eta
_{2}^{0})}}\label{eq:inv3}\\
\eta_{2}^{0}e^{-ta_{2}(\xi_{1}^{0}\eta_{1}^{0},\xi_{2}^{0}\eta_{2}^{0})}  &
=i\left(  \overline{\xi_{2}^{0}e^{ta_{2}(\xi_{1}^{0}\eta_{1}^{0},\xi_{2}%
^{0}\eta_{2}^{0})}}\right)  \label{eq:inv4}%
\end{align}
if we choose $t=0$ then from the above we can see that $\xi_{1}^{0}$ and
$\eta_{1}^{0}$ are real and that $\xi_{2}^{0}=i\overline{\eta_{2}^{0}}.$ Using
the fact that $\xi_{1}^{0}$,$\eta_{1}^{0}\in\mathbb{R}$ with (\ref{eq:inv1})
or (\ref{eq:inv3}) we can see that $a_{1}$ must be real. Using (\ref{eq:inv2})
or (\ref{eq:inv4}) and the fact that $\xi_{2}^{0}=i\overline{\eta_{2}^{0}}$ we
can see that $a_{2}$ is pure imaginary.

All periodic solutions have the initial conditions $\xi_{1}^{0}=\eta_{1}%
^{0}=0.$ If in addition we choose $\xi_{2}^{0}$ of the form $\xi_{2}%
^{0}=re^{i\varphi}$ then for the solution to be real, from (\ref{eq:inv2}), we
must have $\xi_{2}^{0}=i\overline{\eta_{2}^{0}}.$ In such case equation
(\ref{eq: xi_k=...}) gives us the periodic solutions as
\begin{align}
\xi_{2}(t)  &  =\xi_{2}^{0}e^{ta_{2}(0,\xi_{2}^{0}\eta_{2}^{0})}%
=re^{ta_{2}(0,ir^{2})+i\varphi},\label{eq:l(c)(t)-in-xi-eta}\\
\eta_{2}(t)  &  =\eta_{2}^{0}e^{-ta_{2}(0,\xi_{2}^{0}\eta_{2}^{0}%
)}=ire^{-ta_{2}(0,ir^{2})-i\varphi}.\nonumber
\end{align}

\end{proof}

Lemma \ref{lem: a2 imaginary} shows that all periodic solutions of
(\ref{eq: x'=JGradH}) which are real and lie close to the equilibrium point,
are given by the equation
\begin{equation}
l(r,t)=(0,re^{ta_{2}(0,ir^{2})+i\varphi},0,ire^{-ta_{2}(0,ir^{2})-i\varphi}),
\label{eq: lr(t) =}%
\end{equation}
when seen in the $\xi_{\nu},\eta_{\nu}$ coordinates. Let us denote the set
which contains these orbits by
\begin{equation}
B_{R}=\{(0,re^{i\varphi},0,ire^{-i\varphi})|\varphi\in\lbrack0,2\pi),0\leq
r\leq R\}, \label{eq: BR=}%
\end{equation}
where $R$ is sufficiently small for the series $a_{2}(0,ir^{2})$ to be
convergent for $r\leq R.$

Let $P:\mathbb{R}^{4}\rightarrow\mathbb{R}^{4}$ be the time $2\pi$ shift along
the trajectory of (\ref{eq: x'=JGradH}) i.e.
\begin{equation}
P(q(t))=q(t+2\pi), \label{eq:Poincare-shift}%
\end{equation}
where $q(t)$ is a solution of (\ref{eq: x'=JGradH}).

\begin{lemma}
\label{lem: if a_2 ne 0 then P is a twist} If in the series $a_{2}$ from
Theorem \ref{th: Twierdzenie Mosera} i.e.
\begin{equation}
a_{2}(\xi_{1}\eta_{1},\xi_{2}\eta_{2})=\alpha_{2}+a_{2,1}\xi_{1}\eta
_{1}+a_{2,2}\xi_{2}\eta_{2}+... \label{eq:a2-series}%
\end{equation}
for the coefficient $a_{2,2}\in\mathbb{R}$ we have $a_{2,2}\neq0,$ then for a
sufficiently small $R,$ the time $2\pi$ shift along the trajectory $P$
restricted to the set $B_{R}$ is an analytic twist map i.e.
\begin{align}
P(r,\varphi)  &  =(r,\varphi+f(r))\label{eq:Poincare-(r,theta)}\\
\frac{df}{dr}  &  \neq0.\nonumber
\end{align}

\end{lemma}

\begin{proof}
In the $\xi,\eta$ coordinates on $B_{R}$ from (\ref{eq: lr(t) =}) we can see
that the map $P$ takes form
\begin{equation}
P\left(  0,re^{i\varphi},0,ire^{-i\varphi}\right)  =\left(  0,re^{2\pi
a_{2}(0,ir^{2})+i\varphi},0,ire^{-2\pi a_{2}(0,ir^{2})-i\varphi}\right)
.\nonumber
\end{equation}
Keeping in mind that $a_{2}$ is pure imaginary we can see that $P(r,\varphi
)=(r,\varphi-i2\pi a_{2}(0,ir^{2}))$. Since $a_{2}(0,ir^{2})=\alpha
_{2}+a_{2,2}ir^{2}+O(r^{4})$ it is evident that if $a_{2,2}\neq0,$ then for
sufficiently small $r$%
\begin{equation}
\frac{d}{dr}\left(  a_{2}(0,ir^{2})\right)  \neq0. \label{eq:da2-not-zero}%
\end{equation}

\end{proof}

\begin{remark}
\label{rem:twist-action-angle}Note that from Lemma
\ref{lem: if a_2 ne 0 then P is a twist} follows the twist property in
\textit{action--angle coordinates }$(I,\varphi)$ with $I=r^{2}/2.$
\textit{This observation will play a role when applying the KAM Theorem
\ref{th:KAM}. }Observe also that in\textit{ action-angle coordinates the map
}$P$\textit{ has the following form }$P(I,\varphi)=\left(  I,\varphi
+2\pi\mathrm{im}(\alpha_{2})+a_{2,2}I+O(I^{2})\right)  $\textit{, which means
that and the twist is more uniform and does not converge to zero as}
$I\rightarrow0$\textit{, because }$\frac{\partial P_{\varphi}}{\partial
I}(I,\varphi)\rightarrow2\pi a_{22}$\textit{ fo}r $I\rightarrow0$.
\end{remark}

We will now show how to determine whether for a given problem
(\ref{eq: x'=JGradH}) we have $a_{2,2}\neq0.$ For this we will quickly outline
the construction of Moser \cite{Moser} in order to obtain a formula for
$a_{2,2}$. The construction is performed in the following two steps
\begin{gather*}%
\begin{array}
[c]{ccccc}%
\mathbb{C}^{4} & \overset{\Psi}{\rightarrow} & \mathbb{C}^{4} & \overset{\Phi
}{\rightarrow} & \mathbb{C}^{4},
\end{array}
\\%
\begin{array}
[c]{ccccc}%
(\xi_{1},\xi_{2},\eta_{1},\eta_{2}) & \overset{\Psi}{\mapsto} & (x_{1}%
,x_{2},y_{1},y_{2}) & \overset{\Phi}{\mapsto} & (q_{1},q_{2},p_{1},p_{2}),
\end{array}
\end{gather*}
where the transformation $\Phi$ changes the system (\ref{eq: x'=JGradH}) in
the $q_{1},q_{2},p_{1},p_{2}$ coordinates into a system with a simplified
form
\begin{equation}%
\begin{array}
[c]{l}%
\dot{x}_{\nu}=\alpha_{\nu}x_{\nu}+f_{\nu}(x,y)\\
\dot{y}_{\nu}=-\alpha_{\nu}y_{\nu}+g_{\nu}(x,y)
\end{array}
\quad\nu=1,2, \label{eq: x'=alpha x +f}%
\end{equation}
where $\alpha_{1}$ and $\alpha_{2}$ are the eigenvalues of the equilibrium
point and $f$ and $g$ are power series starting from quadratic terms. From the
simplified form (\ref{eq: x'=alpha x +f}) the transformation $\Psi$ determines
the series from Theorem \ref{th: Twierdzenie Mosera}.

The transformation $\Phi$ is a linear function which changes the coordinates
so that the linear part of the equations (\ref{eq: x'=JGradH}) in the new
coordinates becomes generated by a diagonal matrix. Moreover, the
transformation $\Phi$ should be canonical i.e.
\begin{equation}
\Phi^{T}J\Phi=J \label{eq:Phi-kan}%
\end{equation}
where
\begin{equation}
J=\left(
\begin{array}
[c]{cc}%
0 & Id\\
-Id & 0
\end{array}
\right)  \quad\text{and\quad}Id=\left(
\begin{array}
[c]{cc}%
1 & 0\\
0 & 1
\end{array}
\right)  .\nonumber
\end{equation}
Moreover, $\Phi$ should satisfy the following reality condition
\cite{Siegel,Moser}, which expresses the fact that $J_{w}$ is
simply the map describing how the complex conjugation works in new
coordinates
\begin{equation}
J_{z}\Phi=\Phi J_{w}, \label{eq: reality condition}%
\end{equation}
where
\begin{align}
J_{z}(q_{1},q_{2},p_{1},p_{2})  &  =(\bar{q}_{1},\bar{q}_{2},\bar{p}_{1}%
,\bar{p}_{2})\label{eq:reality-corr}\\
J_{w}(x_{1},x_{2},y_{1},y_{2})  &  =(\bar{x}_{1},i\bar{y}_{2},\bar{y}%
_{1},i\bar{x}_{2}).\nonumber
\end{align}

The construction of the transformation $\Psi$ is done by comparison of
coefficients. We look for
\begin{equation}
\Psi=(\phi_{1}(\xi,\eta),\phi_{2}(\xi,\eta),\psi_{1}(\xi,\eta),\psi_{2}%
(\xi,\eta)),\nonumber
\end{equation}
with power series $\phi_{\nu},\psi_{\nu},a_{v},$ $\nu=1,2$ of the form
\begin{equation}%
\begin{array}
[c]{c}%
\phi_{\nu}(\xi_{1},\xi_{2},\eta_{1},\eta_{2})=\sum_{k=1}^{2}\delta_{\nu k}%
\xi_{k}+h.o.t.\\
\psi_{\nu}(\xi_{1},\xi_{2},\eta_{1},\eta_{2})=\sum_{k=1}^{2}\delta_{\nu k}%
\eta_{k}+h.o.t.
\end{array}
\label{eq: phi_v=... , psi_v=...}%
\end{equation}
such that
\begin{align}
x_{\nu}  &  =\phi_{\nu}(\xi_{1},\xi_{2},\eta_{1},\eta_{2}%
)\label{eq:cond-coeff1}\\
y_{\nu}  &  =\psi_{\nu}(\xi_{1},\xi_{2},\eta_{1},\eta_{2}),\nonumber
\end{align}
satisfy (\ref{eq: x'=alpha x +f}) if
\begin{align}
\dot{\xi}_{k}  &  =a_{k}(\xi_{1}\eta_{1},\xi_{2}\eta_{2})\xi_{k}%
\label{eq:cond-coeff2}\\
\dot{\eta}_{k}  &  =-a_{k}(\xi_{1}\eta_{1},\xi_{2}\eta_{2})\eta_{k}.\nonumber
\end{align}
To construct $\phi_{\nu},\psi_{\nu},a_{\nu}$ one can rewrite using
(\ref{eq:cond-coeff1}) and (\ref{eq:cond-coeff2}) the equation
(\ref{eq: x'=alpha x +f}) as
\begin{equation}%
\begin{array}
[c]{l}%
\dot{x}_{\nu}=\sum_{k=1}^{2}\left(  \frac{\partial\phi_{\nu}}{\partial\xi_{k}%
}a_{k}\xi_{k}-\frac{\partial\phi_{\nu}}{\partial\eta_{k}}a_{k}\eta_{k}\right)
=\alpha_{\nu}\phi_{\nu}+f_{\nu}(\phi,\psi)\\
\dot{y}_{\nu}=\sum_{k=1}^{2}\left(  \frac{\partial\psi_{\nu}}{\partial\xi_{k}%
}a_{k}\xi_{k}-\frac{\partial\psi_{\nu}}{\partial\eta_{k}}a_{k}\eta_{k}\right)
=-\alpha_{\nu}\psi_{\nu}+g_{\nu}(\phi,\psi),
\end{array}
\quad\nu=1,2. \label{eq: comparison of coefficients}%
\end{equation}
and compare the coefficients in (\ref{eq: comparison of coefficients}). Let us
denote by $\phi_{\nu,N},\psi_{\nu,N},a_{\nu,N}$ the coefficients in the series
$\phi_{\nu},\psi_{\nu},a_{\nu}$ which come from the homogenous polynomials of
order $N$. We can rewrite the part of (\ref{eq: comparison of coefficients})
which contains all the terms of order $N$ as
\begin{equation}%
\begin{array}
[c]{l}%
\sum_{k=1}^{2}\alpha_{k}\left(  \xi_{k}\frac{\partial}{\partial\xi_{k}}%
-\eta_{k}\frac{\partial}{\partial\eta_{k}}\right)  \phi_{\nu,N}+\ldots
+\delta_{\nu k}\xi_{k}a_{k,N-1}=\alpha_{\nu}\phi_{\nu,N}+\ldots\\
\sum_{k=1}^{2}\alpha_{k}\left(  \xi_{k}\frac{\partial}{\partial\xi_{k}}%
-\eta_{k}\frac{\partial}{\partial\eta_{k}}\right)  \psi_{\nu,N}+\ldots
-\delta_{\nu k}\eta_{k}a_{k,N-1}=-\alpha_{\nu}\psi_{\nu,N}+\ldots
\end{array}
\label{eq: comparison of coefficients with N}%
\end{equation}
where the dots indicate all the terms which can be computed from $\phi_{\nu
,l},\psi_{\nu,l},a_{\nu,l-1}$ with $l=1,\ldots,N-1.$

The nature of equations (\ref{eq: comparison of coefficients with N}) suggests
that the series can be constructed by induction starting with the lowest
terms. It turns out though that not all of the coefficients can be computed
from (\ref{eq: comparison of coefficients with N}). This is because some of
the terms in (\ref{eq: comparison of coefficients with N}) cancel each other
out. If we consider a homogenous polynomial $c\xi_{1}^{n_{1}}\eta_{1}^{m_{1}%
}\xi_{2}^{n_{2}}\eta_{2}^{m_{2}}$ of order $N$ from $\phi_{\nu,N},$ such term
will cancel out in (\ref{eq: comparison of coefficients with N}) if
\begin{equation}
\sum_{k=1}^{2}\alpha_{k}\left(  \xi_{k}\frac{\partial}{\partial\xi_{k}}%
-\eta_{k}\frac{\partial}{\partial\eta_{k}}\right)  c\xi_{1}^{n_{1}}\eta
_{1}^{m_{1}}\xi_{2}^{n_{2}}\eta_{2}^{m_{2}}-\alpha_{\nu}c\xi_{1}^{n_{1}}%
\eta_{1}^{m_{1}}\xi_{2}^{n_{2}}\eta_{2}^{m_{2}}=0.\nonumber
\end{equation}
This can happen only if we have
\begin{equation}
\sum_{k=1}^{2}\alpha_{k}\left(  n_{k}-m_{k}\right)  -\alpha_{\nu}=0.
\label{eq:poly-coeff=0}%
\end{equation}
By the assumption of the Theorem \ref{th: Twierdzenie Mosera} that for any
$t\in\mathbb{R}$ we have $t\alpha_{1}+\alpha_{2}\neq0,$ we can see that
(\ref{eq:poly-coeff=0}) is true only for the terms of the form $c\xi
_{v}\left(  \xi_{1}\eta_{1}\right)  ^{n_{1}}\left(  \xi_{2}\eta_{2}\right)
^{n_{2}}.$ The value of the coefficient $c$ corresponding to such terms is
chosen from an appropriate normalization condition \cite{Moser}. In our case
though the choice of the normalization does not play an important role. We are
interested in computation of the term $a_{2,2}$ and this can be done by
induction starting with $N=1$ and stopping at $N=3.$ For $N=1$ all the terms
are uniquely determined. For $N=2$ there are no terms which would cancel out
in (\ref{eq: comparison of coefficients with N}). For $N=3$ we can use the
first equation from (\ref{eq: comparison of coefficients with N}) to find
$a_{2,2}.$ There the coefficient $a_{2,2}$ stands together with $\xi
_{2}\left(  \xi_{2}\eta_{2}\right)  $ and the term for $\xi_{2}\left(  \xi
_{2}\eta_{2}\right)  $ in $\phi_{\nu,3}$ will cancel out from the equation.
Therefore $a_{2,2}$ is uniquely determined since it depends only on $\phi
_{\nu,1},\psi_{\nu,1},\phi_{\nu,2}$ and $\psi_{\nu,2}$. Using this procedure
we can obtain a direct formula for the term $a_{2,2}.$

\begin{lemma}
\label{lem: a2,2=...}If
\begin{equation}%
\begin{array}
[c]{c}%
f_{\nu}(x_{1},x_{2},y_{1},y_{2})=\sum_{i,j,k,l\geq1}f_{ijkl}^{\nu}x_{1}%
^{i}x_{2}^{j}y_{1}^{k}y_{2}^{l}\\
g_{\nu}(x_{1},x_{2},y_{1},y_{2})=\sum_{i,j,k,l\geq1}g_{ijkl}^{\nu}x_{1}%
^{i}x_{2}^{j}y_{1}^{k}y_{2}^{l}%
\end{array}
\quad\nu=1,2.\nonumber
\end{equation}
then
\begin{align}
a_{2,2}  &  =\frac{1}{\alpha_{2}}(-f_{1,1,0,0}^{2}f_{0,1,0,1}^{1}%
-f_{0,1,1,0}^{2}g_{0,1,0,1}^{1}+f_{0,0,1,1}^{2}g_{0,2,0,0}^{1}-f_{0,2,0,0}%
^{2}f_{0,1,0,1}^{2}\label{eq: formula for a22}\\
&  \quad+2g_{0,2,0,0}^{2}f_{0,0,0,2}^{2}+f_{1,0,0,1}^{2}f_{0,2,0,0}%
^{1}-g_{0,1,0,1}^{2}f_{0,1,0,1}^{2})+f_{0,2,0,1}^{2}\nonumber
\end{align}

\end{lemma}

\begin{proof}
The above can be checked from the formula
(\ref{eq: comparison of coefficients with N}) by direct computation.
\end{proof}

This ends our construction of the coefficient $a_{2,2}$.

Let us now briefly turn to the relation between the energy and the radius $r$
of the periodic orbits (\ref{eq: lr(t) =}).

\begin{lemma}
\label{lem:dist(Lc,L2)}For sufficiently small $r$ the energy of the orbit
$l_{r}$ (\ref{eq: lr(t) =}) i.e.
\begin{equation}
h(r)=H(\Phi\left(  \Psi(l_{r})\right)  ), \label{eq:h(r)=C2+rkwadrat...}%
\end{equation}
is equal to
\begin{equation}
h(r)=H(0)+\frac{1}{2}D^{2}H(0)\left(  \Phi(0,1,0,i)\right)  r^{2}+o(r^{2}),
\label{eq:h(r)=...h2}%
\end{equation}
where $D^{2}H(0)\left(  \Phi(0,1,0,i)\right)  $ denote the value of the
quadratic form $D^{2}H(0)$ on the vector $\Phi(0,1,0,i)$.
\end{lemma}

\begin{proof}
Since the problem (\ref{eq: x'=JGradH}) is autonomous the energy is constant
along the orbit $l_{r}.$ Without any loss of generality we can therefore
assume that $\varphi$ in equation (\ref{eq: lr(t) =}) for $l_{r}$ is zero and
compute
\begin{equation}
h(r)=H(\Phi\left(  \Psi(l(r,0))\right)  ).\nonumber
\end{equation}
Let us first note that the construction of $\Psi=(\phi_{1},\phi_{2},\psi
_{1},\psi_{2})$ produced power series of the form
(\ref{eq: phi_v=... , psi_v=...}), hence
\begin{equation}
\Psi(l(r,0))=(\phi_{1},\phi_{2},\psi_{1},\psi_{2})(l(r,0))=(0,r,0,ir)+O(r^{2}%
).\nonumber
\end{equation}
The transformation $\Phi$ is linear and therefore
\begin{equation}
\Phi\left(  \Psi(l(r,0))\right)  =r\Phi(0,1,0,i)+O(r^{2}).
\label{eq:Psi-radius}%
\end{equation}
We can compute $h(r)$ as
\begin{align}
h(r)  &  =H(\Phi\left(  \Psi(l(r,0))\right)  ) \label{eq:proof-h(r)-finishing}%
\\
&  =H(0)+DH(0)\left(  \Phi\left(  \Psi(l(r,0))\right)  \right) \nonumber\\
&  \quad+\frac{1}{2}D^{2}H(0)\left(  \Phi\left(  \Psi(l(r,0))\right)  \right)
^{2}+o(\left\vert \Phi\left(  \Psi(l(r,0))\right)  \right\vert ^{2}).\nonumber
\end{align}
Since zero is an equilibrium point we know that $DH(0)=0,$ thus by
substituting (\ref{eq:Psi-radius}) into (\ref{eq:proof-h(r)-finishing}) we
obtain our claim.
\end{proof}

\section{Twist in the PRC3BP at $L_{2}^{\mu}$}

\label{sec:twist-L2}

As mentioned in Section~\ref{sec:prelim}, we have a family of periodic
Lyapounov orbits $l_{\mu}(C)$ around $L_{2}^{\mu}$ for energies larger than
and sufficiently close to the energy $C_{2}^{\mu}$ of the libration point
$L_{2}^{\mu}.$ This family of orbits corresponds to the set $B_{R}$ (see
(\ref{eq: BR=})) of orbits constructed in the previous section. In this
section we will show that for sufficiently small $\mu$ we have a twist
property on this family of periodic orbits. The main idea for the proof is to
approximate the PRC3BP (with Hamiltonian (\ref{eq:PRC3BP-Hill})) expressed in
Hill' coordinates (\ref{eq:Hill-coordinates}) with the Hill's problem
(\ref{eq: Hills problem}). First we shall consider the Hill's problem
(\ref{eq: Hills problem}) where we have explicit formulas for the libration
point $L_{2}^{\text{Hill}}$, the linearized vector field at $L_{2}%
^{\text{Hill}}$ its eigenvalues etc., which will allow us to compute the twist
coefficient $a_{2,2}^{\text{Hill}}$. We will then show that the coefficients
$a_{2,2}^{\mu}$ computed for the PRC3BP in Hill's coordinates, converges to
$a_{2,2}^{\text{Hill}}$ as $\mu$ tends to zero.

\begin{lemma}
\label{lem:twist-hill}Let $P^{\text{Hill}}$ be the time $2\pi$ shift along the
trajectory Poincar\'e map of the Hill's problem (\ref{eq: Hills problem}).
Then there exists a radius $R_{\mathrm{Hill}}\in\mathbb{R}$, such that the map
$P^{\text{Hill}}$ expressed in in radius angle coordinates on the set of
Lyapounov orbits around $L_{2}^{\mathrm{Hill}}$
\begin{equation}
P^{\text{Hill}}:\left[  0,R_{\mathrm{Hill}}\right]  \times S^{1}%
\rightarrow\left[  0,R_{\mathrm{Hill}}\right]  \times S^{1},\nonumber
\end{equation}
is a twist map.
\end{lemma}

\begin{proof}
We will apply the procedure from the previous section and compute
$a_{2,2}^{\text{Hill}}$ for the equilibrium point $L_{2}^{\text{Hill}%
}=(3^{-1/3},0,0,3^{-1/3}).$ The linear terms of (\ref{eq: Hills problem}) in
$L_{2}^{\text{Hill}}$ are given by (\ref{eq:Hill-linear}) with eigenvalues
$\pm\alpha_{1}$, $\pm\alpha_{2}$, $\alpha_{1}=\sqrt{1+2\sqrt{7}}$ and
$\alpha_{2}=\sqrt{1-2\sqrt{7}}$. The first of the two is real and the second
is pure imaginary. We will choose the function $\Phi^{\text{Hill}}$ composed
of the eigenvectors of the eigenvalues $\pm\alpha_{1}$ and $\pm\alpha_{2}$%
\begin{equation}
\Phi^{\text{Hill}}=\left(
\begin{array}
[c]{cccc}%
\lambda_{1} & \beta & \lambda_{2} & -i\bar{\beta}\\
-9\lambda_{1}\frac{1}{\alpha_{1}\left(  \sqrt{7}+4\right)  } & -\beta9\frac
{1}{\alpha_{2}\left(  \sqrt{7}-4\right)  } & 9\lambda_{2}\frac{1}{\alpha
_{1}\left(  \sqrt{7}+4\right)  } & -i\bar{\beta}9\frac{1}{\alpha_{2}\left(
\sqrt{7}-4\right)  }\\
9\lambda_{1}\frac{\sqrt{7}+3}{\alpha_{1}\left(  \sqrt{7}+4\right)  } &
-\beta9\frac{\sqrt{7}-3}{\alpha_{2}\left(  \sqrt{7}-4\right)  } &
-9\lambda_{2}\frac{\sqrt{7}+3}{\alpha_{1}\left(  \sqrt{7}+4\right)  } &
-i\bar{\beta}9\frac{\sqrt{7}-3}{\alpha_{2}\left(  \sqrt{7}-4\right)  }\\
-\lambda_{1}\frac{2}{3+\sqrt{7}} & \beta\frac{2}{\sqrt{7}-3} & -\lambda
_{2}\frac{2}{3+\sqrt{7}} & -i\bar{\beta}\frac{2}{\sqrt{7}-3}%
\end{array}
\right)  , \label{eq: Phi=}%
\end{equation}
with $\lambda_{1},\lambda_{2}\in\mathbb{R},$ $\beta\in\mathbb{C}$.
The above transformation $\Phi^{\text{Hill}}$ satisfies the
reality condition (\ref{eq: reality condition}), and if  we choose
the coefficients $\lambda_{1},\lambda_{2},\beta$ as
\begin{align}
\lambda_{1}  &  =-\lambda_{2}=\frac{1}{6}\sqrt{\frac{\alpha_{1}\left(
\sqrt{7}+4\right)  }{\sqrt{7}}},\label{eq:lambda-beta}\\
\beta &  =\frac{1}{6}\sqrt{\frac{i\alpha_{2}\left(  \sqrt{7}-4\right)  }%
{\sqrt{7}}},\nonumber
\end{align}
then $\Phi^{\text{Hill}}$ is also canonical. Computing the power series
$f_{\nu}$ and $g_{\nu}$ from Lemma \ref{lem: a2,2=...} at $L_{2}^{\text{Hill}%
}$ and using (\ref{eq: formula for a22}) to compute the term $a_{2,2}%
^{\text{Hill}},$ by rather laborious computations (performed in Maple) one
will obtain
\begin{equation}
a_{2,2}^{\text{Hill}}=\frac{\sqrt[3]{9}}{224}\left(  102\sqrt{7}-57\right)
\approx\allowbreak1.\,\allowbreak976\,7.\nonumber
\end{equation}
Since $a_{2,2}^{\text{Hill}}\neq0$, by Lemma
\ref{lem: if a_2 ne 0 then P is a twist} we have the twist property in the
radius angle coordinates for all $r$ such that $0<r<R_{\text{Hill}}$, where
$R_{\text{Hill}}$ is sufficiently small.
\end{proof}

\begin{remark}
Let us stress that $R_{\text{Hill}}$ is independent of $\mu$. From our
attempts of rigorous estimation of $R_{Hill}$, following the estimates
conducted during the proof of Theorem~\ref{th: Twierdzenie Mosera} in
\cite{Moser}, we obtain that that $R_{Hill}\approx10^{-4}$.
\end{remark}

One can also apply the procedure outlined in Section \ref{sec:lap-twist} to
compute the coefficient $a_{2,2}^{\mu}$ for the PRC3BP with Hamiltonian
(\ref{eq:PRC3BP-Hill}) expressed in Hill' coordinates
(\ref{eq:Hill-coordinates}). To do so one needs to compute the libration point
$\bar{L}_{2}^{\mu}=\mu^{-1/3}\left(  L_{2}^{\mu}-\left(  \mu-1,0,0,\mu
-1\right)  \right)  $, compute the expansion of the vector field at $\bar
{L}_{2}^{\mu_{k}}$ up to the order three, compute the linear change of
coordinates $\Phi=\Phi(\mu)$ and compute $a_{2,2}^{\mu} $ using Lemma
\ref{lem: a2,2=...}. Numerical results of such computations for a selection of
parameters from the family $\{\mu_{k} \}$ from Theorem
\ref{lem homoclinic orbit for uk} are given in the below table. The values
$\mu_{k}$ chosen in the table are the numerical approximations of the series
(\ref{eq:muk}) from \cite{Simo}.

\begin{samepage}
\begin{equation}%
\begin{tabular}{r | r | r}%
$k$ & $\mu_{k}$ & $ a_{2,2}^{\mu_k}$\\ \hline
2 & 0.4253863522E-2  & 1.967649155\\
3 & 0.6752539971E-3  & 1.968237635\\
4 & 0.2192936884E-3  & 1.970039000\\
&    & \\
10 & 0.92907436E-5  & 1.973971883\\
11 & 0.68212830E-5  & 1.974219993\\
12 & 0.51549632E-5  & 1.974426964\\
&    & \\
50 & 0.582146E-7  & 1.976164111\\
60 & 0.336890E-7  & 1.976252225\\
70 & 0.212152E-7  & 1.976315175\\
&    & \\
200 & 0.9096E-9  & 1.976563023
\end{tabular}
\notag \end{equation}
\begin{equation}
a_{2,2}^{Hill} \approx 1.9767 \notag
\end{equation}
\begin{center}
\textbf{Table 1. }The twist coefficient $a_{2,2}$ for various
masses $\mu _{k}.$\bigskip
\end{center}
\end{samepage}

\begin{theorem}
\label{lem:twist-c}Let $P^{\mu}$ denote the the time $2\pi$ shift along the
trajectory $P^{\mu}$ of the PRC3BP in Hill's coordinates
(\ref{eq:Hill-coordinates}). Then for sufficiently small $\mu$ the map
$P^{\mu}$ expressed in in radius angle coordinates on the set of Lyapounov
orbits around $\bar{L}_{2}^{\mu} =\mu^{-1/3}\left(  L_{2}^{\mu}-\left(
\mu-1,0,0,\mu-1\right)  \right)  $ is an analytic twist map i.e. for
$r<R_{\text{Hill}}$
\begin{equation}
P^{\mu}(r,\varphi)=(r,\varphi+f(r))\nonumber
\end{equation}
and
\begin{equation}
\frac{df}{dr}\neq0\quad\text{for all }r\in\lbrack0,R_{\text{Hill}}].\nonumber
\end{equation}

\end{theorem}

\begin{proof}
We observe that $\bar{L}_{2}^{\mu}$ depends analytically on $\mu^{1/3}$. First
let us note that since the operator $\Phi$ from our construction brings the
derivative of the vector field to the Jordan form, the operator $\Phi
_{3\text{Body}}^{\mu}$ for the PRC3BP can be chosen close (depending
analytically on $\mu^{1/3}$) to the operator $\Phi^{\text{Hill}}.$ The same
can be said about the coefficients $f_{\nu}$ and $g_{\nu}$ from Lemma
\ref{lem: a2,2=...}, since those come from the Taylor expansion up to the
third order of the vector field at $\bar{L}_{2}^{\mu}.$ Moreover, from the
proof of the Lyapounov Moser Theorem \ref{th: Twierdzenie Mosera} in
\cite{Moser}, we know that the radius of convergence of the series from the
Theorem \ref{th: Twierdzenie Mosera} can be chosen to be independent from
$\mu,$ for $\mu$ sufficiently close to zero. This means that the coefficient
$a_{2,2}^{\mu}$ constructed for he PRC3BP will tend to the coefficient
$a_{2,2}^{\text{Hill}}$ of the Hill's problem
\begin{equation}
\lim_{\mu\rightarrow0}a_{2,2}^{\mu}=a_{2,2}^{\text{Hill}}\thickapprox
1.\,\allowbreak976\,7. \label{eq:a2-limit}%
\end{equation}

\end{proof}

\begin{remark}
From Theorem \ref{lem:twist-c} and the results from Table 1, it is reasonable
to believe that we will have twist for all $\mu_{k}$ for $k\geq2$. Let us
point out that in the Hill's coordinates the radius of convergence can be
chosen to be independent from $\mu.$ In the original coordinates of the PRC3BP
though, since $r_{\text{Hill}}=\mu^{1/3}r,$ this radius will depend and
decrease with $\mu$ as $R^{\mu}=\mu^{1/3}R_{\text{Hill}}$. This means that in
the original coordinates we will have the twist property only for orbits with
a radius smaller than $\mu^{1/3}R_{\text{Hill}}.$
\end{remark}

\section{Normal hyperbolicity, KAM theorem and the persistence of Lyapounov
orbits}

\label{sec:norm-hyp-KAM}

In this Section we briefly recall some facts from the normal hyperbolicity
theory and a version of the K.A.M. (Kolmogorov, Arnold, Moser) Theorem. We
then apply the results to obtain the persistence result of a Cantor family of
Lyapounov orbits around $L_{2}^{\mu}$ for the perturbation from PRC3BP to
PRE3BP. Our approach closely follows the method of \cite{Llave}. We will
therefore rewrite the theorems used in \cite{Llave} and verify that their
assumptions are satisfied in our particular setting.

First let us recall the results concerning normal hyperbolicity.

\begin{definition}
[{\cite[A1]{Llave}}]\label{def:norm-hyp} Let $M$ be a manifold in
$\mathbb{R}^{n}$ and $\Phi_{t}$ a $C^{r},$ $r\geq1$ flow on it. We say that a
(smooth) manifold $\Lambda\subset M$ -- possibly with boundary -- invariant
under $\Phi_{t}$ is $\alpha$-$\beta$ normally hyperbolic when there is a
bundle decomposition
\begin{equation}
TM=T\Lambda\oplus E^{s}\oplus E^{u},\nonumber
\end{equation}
invariant under the flow, and numbers $C>0,$ $0<\beta<\alpha,$ such that for
$x\in\Lambda$%
\begin{align}
v  &  \in E_{x}^{s}\Leftrightarrow|D\Phi_{t}(x)v|\leq Ce^{-\alpha t}%
|v|\quad\forall t>0,\label{eq:norm-hyp-cond1}\\
v  &  \in E_{x}^{u}\Leftrightarrow|D\Phi_{t}(x)v|\leq Ce^{\alpha t}%
|v|\quad\forall t<0,\label{eq:norm-hyp-cond2}\\
v  &  \in T_{x}\Lambda\Leftrightarrow|D\Phi_{t}(x)v|\leq Ce^{\beta|t|}%
|v|\quad\forall t. \label{eq:norm-hyp-cond3}%
\end{align}

\end{definition}

\begin{theorem}
[{\cite[A7]{Llave}}]\label{th:norm-hyp-manifolds-exist}Let $\Lambda$ be a
compact $\alpha$-$\beta$ normally hyperbolic manifold (possibly with a
boundary) for the $C^{r}$ flow $\Phi_{t},$ satisfying the Definition
\ref{def:norm-hyp}. Then there exists a sufficiently small neighborhood
$U\ $of $\Lambda$ and a sufficiently small $\delta>0$ such that

\begin{enumerate}
\item The manifold $\Lambda$ is $C^{\text{min}(r,r_{1}-\delta)},$ where
$r_{1}=\alpha/\beta.$

\item For any $x$ in $\Lambda,$ the set
\begin{align*}
W_{x}^{s}  &  =\{y\in U:\text{dist}(\Phi_{t}(y),\Phi_{t}(x))\leq
Ce^{(-\alpha+\delta)t}\text{ for }t>0\}\\
&  =\{y\in U:\text{dist}(\Phi_{t}(y),\Phi_{t}(x))\leq Ce^{(-\beta-\delta
)t}\text{ for }t>0\}
\end{align*}
is a $C^{r}$ manifold and $T_{x}W_{x}^{s}=E_{x}^{s}.$

\item The bundles $E_{x}^{s}$ are $C^{\text{min}(r,r_{0}-\delta)}$ in $x,$
where $r_{0}=(\alpha-\beta)/\beta,$ and
\begin{align*}
W_{\Lambda}^{s}  &  =\{y\in U:\text{dist}(\Phi_{t}(y),\Lambda)\leq
Ce^{(-\alpha+\delta)t}\text{ for }t>0\}\\
&  =\{y\in U:\text{dist}(\Phi_{t}(y),\Lambda)\leq Ce^{(-\beta-\delta)t}\text{
for }t>0\}
\end{align*}
is a $C^{\min(r,r_{0}-\delta)}$ manifold. Moreover $T_{x}W_{\Lambda}^{s}%
=E_{x}^{s}.$ Finally
\begin{equation}
W_{\Lambda}^{s}=\bigcup_{x\in\Lambda}W_{x}^{s}.\nonumber
\end{equation}
Moreover, we can find a $\rho>0$ sufficiently small and a $C^{\min
(r,r_{0}-\delta)}$ diffeomorphism from the bundle of balls of radius $\rho$ in
$E_{\Lambda}^{s}$ to $W_{\Lambda}^{s}\cap U.$
\end{enumerate}
\end{theorem}

\begin{remark}
An analogous theorem can be stated for $W_{\Lambda}^{u}$ by considering the
flow $\Phi_{-t}.$
\end{remark}

The following Theorem and two Remarks concern the persistence of the normally
hyperbolic manifold and its stable and unstable manifolds.

\begin{theorem}
[{\cite[A.14]{Llave}}]\label{th:pert-norm-hyp} Let $\Lambda\subset M$
($\Lambda$ not necessarily compact) be $\alpha$-$\beta$ normally hyperbolic
for the flow $\Phi_{t}$ generated by the vector field $X,$ which is uniformly
$C^{r}$ in a neighborhood $U$ of $\Lambda$ such that dist$(M\setminus
U,\Lambda)>0.$ Let $\Psi_{t}$ be the flow generated by another vector field
$Y$ which is $C^{r}$ and sufficiently $C^{1}$ close to $X.$ Then we can find a
manifold $\Gamma$ which is $\alpha^{\prime}$-$\beta^{\prime}$ normally
hyperbolic for $Y$ and $C^{\min(r,r_{1}-\delta)}$ close to $\Lambda,$ where
$r_{1}=\alpha/\beta.$

The constants $\alpha^{\prime}$, $\beta^{\prime}$ are arbitrarily close to
$\alpha$, $\beta$ if $Y$ is sufficiently $C^{1}$ close to $X.$

The manifold $\Gamma$ is the only $C^{\min(r,r_{1}-\delta)}$ normally
hyperbolic manifold $C^{0}$ close to $\Lambda$ and locally invariant under the
flow of $Y.$
\end{theorem}

The above Theorem is extended to give us a smooth dependence on the parameter
by the following two remarks.

\begin{remark}
[{see \cite[observation 1. page 390]{Llave}}]\label{rem:pert-norm-hyp} Assume
that we have a family of flows $\Phi_{t,e},$ generated by vector fields
$X_{e}$ which are jointly $C^{r}$ in all its variables (the base point $x$ and
the parameter $e$). Let $\Lambda_{e}$ be the normally hyperbolic manifold
$\Gamma$ from Theorem \ref{th:pert-norm-hyp} for the flow $\Phi_{t,e}.$ Then
there exists a $C^{\min(r,r_{1}-\delta)}$ mapping $F:\Lambda\times
I\rightarrow M,$ where $r_{1}=\alpha/\beta$ and $I\subset\mathbb{R}$ is an
interval containing zero, such that $F(\Lambda,e)=\Lambda_{e}$ and
$F(\cdot,0)$ is the identity.
\end{remark}

\begin{remark}
[{see \cite[observation 2. page 390]{Llave}}]%
\label{rem:pert-norm-hyp-with-manifolds} For a family of flows $\Phi_{t,e}$
with the same assumptions as in Remark \ref{rem:pert-norm-hyp}, there exists a
$C^{\min(r,r_{1}-\delta)}$ ( $r_{1}=\alpha/\beta$ ) mapping $R^{s}:W_{\Lambda
}^{s}\times I\rightarrow M$ such that $R^{s}(W_{\Lambda}^{s},e)=W_{\Lambda
,e}^{s}$, $R^{s}(\cdot,e)|_{\Lambda}=F(\cdot,e),$ $R^{s}(W_{x}^{s}%
,e)=W_{F(x,e),e}^{s}.$

An analogous mapping $R^{u}$ also exists for $W_{\Lambda}^{u}.$
\end{remark}

Let us now turn to a quantitative version of the KAM Theorem used in
\cite{Llave}. Let us recall that a real number $\omega$ is called a
Diophantine number of exponent $\theta$ if there exists a constant $C>0$ such
that
\begin{equation}
\left\vert \omega-\frac{p}{q}\right\vert \geq\frac{C}{q^{\theta+1}}\nonumber
\end{equation}
for all $p\in\mathbb{Z},$ $q\in\mathbb{N}$.

\begin{definition}
Let $(M,\omega_{M}),$ $(N,\omega_{N})$ be two symplectic manifolds of same
dimension. If $(M,\omega_{M}),$ $(N,\omega_{N})$ are exact symplectic (i.e.
there exist one-forms $\alpha_{M},\alpha_{N}$ such that $\omega_{M}%
=d\alpha_{M},$ $\omega_{N}=d\alpha_{N}$) then we say that a diffeomorphism%
\begin{equation}
F:M\rightarrow N\nonumber
\end{equation}
is exact symplectic when there exists a real valued function $G$ on $M$ such
that%
\begin{equation}
F^{\ast}\alpha_{N}-\alpha_{M}=dG.\nonumber
\end{equation}

\end{definition}

\begin{theorem}
[{KAM Theorem \cite[Theorem 4.8]{Llave}}]\label{th:KAM}Let $f:\left[
0,1\right]  \times\mathbb{T}^{1}\rightarrow\left[  0,1\right]  \times
\mathbb{T}^{1}$ be an exact symplectic $C^{l}$ map with $l\geq6.$

Assume that $f=f_{0}+ef_{1}$, where $e\in\mathbb{R,}$
\begin{equation}
f_{0}\left(  I,\varphi\right)  =\left(  I,\varphi+A\left(  I\right)  \right)
, \label{eq:KAM-form}%
\end{equation}
$A$ is $C^{l},$ $\left\vert \frac{dA}{dI}\right\vert \geq M$, and $\left\Vert
f_{1}\right\Vert _{C^{l}}\leq1.$

Then, if $e^{1/2}M^{-1}=\rho$ is sufficiently small, for a set of Diophantine
numbers $\sigma$ of exponent $\theta=5/4$, we can find invariant tori which
are the graph of $C^{l-3}$ functions $u_{\sigma}$, the motion on them is
$C^{l-3}$ conjugate to the rotation by $\sigma$ and the tori cover the whole
annulus except a set of measure smaller than $O\left(  M^{-1}e^{1/2}\right)
.$

Moreover, we can find expansions
\begin{equation}
u_{\sigma}=u_{\sigma}^{0}+eu_{\sigma}^{1}+r_{\sigma}, \label{eq:kam-tor-close}%
\end{equation}
with $u_{\sigma}^{0}=A^{-1}(\sigma)$, $\left\Vert r_{\sigma}\right\Vert
_{C^{l-4}}\leq O\left(  e^{2}\right)  ,$ and $\left\Vert u_{\sigma}%
^{1}\right\Vert _{C^{l-4}}\leq O(1)$.
\end{theorem}

All of the above results have been taken from \cite{Llave}. Now we will apply
them to the setting of the PRE3BP. We will first show that the set of the
Lyapounov orbits of the PRC3BP is normally hyperbolic.

Let $\phi_{t,s}^{e}:\mathbb{R}^{4}\rightarrow\mathbb{R}^{4}$ be given by
\begin{equation}
\phi_{t,s}^{e}(x)=q(s+t),\nonumber
\end{equation}
where $q(\cdot)$ is the solution for the PRE3BP in Hill's coordinates
(\ref{eq:Hill-coordinates}) (generated by the Hamiltonian
(\ref{eq:PRC3BP-Hill})), with an initial condition $q(s)=\mathbf{\bar{x}}$. We
will define the flow on the extended phase space $\Phi_{t}^{e}:\mathbb{R}%
^{4}\times\mathbb{R}\rightarrow\mathbb{R}^{4}\times\mathbb{R}$ as
\begin{equation}
\Phi_{t}^{e}(\mathbf{\bar{x}},s)=(\phi_{t,s}^{e}(\mathbf{\bar{x}}),s+t).
\label{eq:Phi_t_e}%
\end{equation}
Observe that the flow $\Phi_{t}^{e}$ is $2\pi$ periodic with respect to $s$
variable, hence may be equivalently treated as a flow on $\mathbb{R}^{4}\times
S^{1}$. This will later give us uniform $C^{r}$ estimates. We shall use a
notation
\begin{equation}
l(r)=l(r,S^{1})\times S^{1}\nonumber
\end{equation}
to denote a torus of all trajectories of Lyapounov orbits $l(r,t)$ (see
(\ref{eq: lr(t) =})) of radius $r$ in the extended phase space.

\begin{lemma}
\label{lem:Lap-orb normally hyp.}For a sufficiently small mass $\mu$ the set
\begin{equation}
\Lambda=\{l(r)|r\in\lbrack0,R_{\mathrm{Hill}})\}\nonumber
\end{equation}
of Lyapounov orbits of the PRC3BP (in the extended phase space) is $\alpha
$-$\beta$ normally hyperbolic, where $\alpha>0$ is close to the real
eigenvalue $\alpha_{1}$ at the Libration point $\bar{L}_{2}^{\mu}=\mu
^{-1/3}\left(  L_{2}^{\mu}-\left(  \mu-1,0,0,\mu-1\right)  \right)  $ and
$\beta>0$ can be chosen arbitrarily close to zero.
\end{lemma}

\begin{proof}
Consider the $\xi_{1},\xi_{2},\eta_{1},\eta_{2}$ coordinates from the previous
section together with a time $t$ coordinate. The $\xi_{1},\eta_{1}$ are the
coordinates of the hyperbolic expansion and $\xi_{1},\eta_{1}$ are the
coordinates of the twist rotation around the libration point $\bar{L}_{2}%
^{\mu}$. We have
\begin{multline*}
M=\{(\xi_{1},\xi_{2},\eta_{1},\eta_{2},t)|\xi_{1},\eta_{1}\in\mathbb{R}%
,\xi_{2}=re^{i\varphi},\eta_{2}=ire^{-i\varphi},\\
r\in\lbrack0,R_{\mathrm{Hill}}),\varphi\in\lbrack0,2\pi),t\in\lbrack0,2\pi)\}.
\end{multline*}
We can define
\begin{align*}
E^{u}  &  =\{(\xi_{1},0,0,0,0)|\xi_{1}\in\mathbb{R}\},\\
E^{s}  &  =\{(0,0,\eta_{1},0,0)|\eta_{1}\in\mathbb{R}\},\\
T\Lambda &  =\{(0,\xi_{2},0,\eta_{2},t)|\xi_{2}=re^{i\varphi},\eta
_{2}=ire^{-i\varphi},r\in\mathbb{R,}\varphi\in\lbrack0,2\pi),t\in\lbrack
0,2\pi)\},
\end{align*}
then we will have $TM=E^{u}\oplus E^{s}\oplus T\Lambda.$ The conditions
(\ref{eq:norm-hyp-cond1}) and (\ref{eq:norm-hyp-cond2}) are satisfied with a
coefficient $\alpha>0$ close to the eigenvalue $\alpha_{1}$ at $\bar{L}%
_{2}^{\mu}$ because the coordinates $\xi_{1}$ and $\eta_{1}$ are the
coordinates of hyperbolic expansion and contraction. For sufficiently small
$\mu$ the eigenvalue $\alpha_{1}$ is close to $\alpha_{1}^{Hill}%
=\sqrt{1+2\sqrt{7}}\approx\allowbreak2.\,\allowbreak508\,3$ of the Hill's problem.

Let $\Phi_{t}$ be the flow in the extended phase space. From
(\ref{eq: lr(t) =}) for $x=(0,\xi_{2},0,\eta_{2},t)\in\Lambda$ with $\xi
_{2}=re^{i\varphi},$ $\eta_{2}=ire^{-i\varphi}$ we have
\begin{equation}
\left\Vert D\Phi_{t}(x)\right\Vert \leq1+r\left\vert \frac{d}{dr}%
a_{2}(0,ir^{2})\right\vert t,\nonumber
\end{equation}
where $a_{2}$ is the function given by (\ref{eq:a-series}) in Theorem
\ref{th: Twierdzenie Mosera}. The growth of derivative of $D\Phi_{t}(x)$ is at
most linear in $t$. For any $\beta>0$ there exists a constant $C>0$, such that
for all $v\in T_{x}\Lambda$ and all $t$
\begin{equation}
|D\Phi_{t}(x)v|\leq Ce^{\beta|t|}|v|. \label{eq:twistnorm-hyp-cond3}%
\end{equation}

\end{proof}

We will now define the time $2\pi$ shift along a trajectory Poincar\'{e} map
and later apply the KAM Theorem to it. By Lemma
\ref{lem:Lap-orb normally hyp.} for $e=0$ we have an $\alpha$-$\beta$ normally
hyperbolic invariant manifold for $\Phi_{t}^{0}$ of the form $\Lambda
=\{l(r)|r\in\lbrack0,R_{\text{Hill}})\}$. Let $U$ be an open neighborhood of
$\Lambda$. We will define time $2\pi$ shift along a trajectory Poincar\'{e}
map $P_{t_{0}}^{e}: U\cap\{t=t_{0}\}\rightarrow\mathbb{R}^{4}$ as
\begin{equation}
P_{t_{0}}^{e}(x)=\phi_{2\pi,t_{0}}^{e}(x).\nonumber
\end{equation}

We are now ready to apply Theorems \ref{th:pert-norm-hyp} and \ref{th:KAM} to
obtain the following persistence result for the family of the Lyapounov orbits.

\begin{theorem}
\label{th:KAM for L2}Let $R<R_{\mathrm{Hill}}$ be a fixed number and $\kappa$
be the parameter from (\ref{eq:kappa-e-bound}). If we choose sufficiently
small $\mu^{\ast}>0$ then for all $\mu>0$ and $e>0$ for which $e\mu
^{-2/3}<\kappa$ the normally hyperbolic manifold (with a boundary; considered
in the extended phase space) $\Lambda=\{l(r)|r\in\lbrack0,R)\}$ of the PRC3BP
persists under the perturbation to PRE3BP with the parameter $e$, to a
normally hyperbolic manifold (with a boundary) $\Lambda_{e}$. Moreover, for
any such $e$ there exists a Cantor set $\mathcal{C}\subset\lbrack0,R]$ such
that for any $r\in\mathcal{C}$ there exists an invariant (two dimensional)
torus
\begin{equation}
l^{e}(r)=\left\{  l_{t_{0}}^{e}\left(  r\right)  |t_{0}\in S^{1}\right\}
,\nonumber
\end{equation}
where $l_{t_{0}}^{e}\left(  r\right)  $ is an one--dimensional torus invariant
under the map $P_{t_{0}}^{e}$. The family of tori $l^{e}(r)$ for
$r\in\mathcal{C}$, covers $\Lambda_{e}$ except a set of a measure smaller than
$O(\left(  e\mu^{-1/3}\right)  ^{1/2})$.
\end{theorem}

\begin{proof}
The PRE3BP in Hill's coordinates (\ref{eq:Hill-coordinates}) is generated by
the Hamiltonian $\bar{H}^{e}$ from (\ref{eq:He-reduced}). In the neighborhood
$U$ of $\Lambda$ the $\bar{r} $ from (\ref{eq:PRC3BP-Hill-reduced}) and
(\ref{eq:He-reduced}) is bounded and separated from zero, which means that the PRE3BP is a uniform
perturbation of the Hill's problem (\ref{eq: Hills problem}). By Lemma
\ref{lem:Lap-orb normally hyp.} we know that $\Lambda$ is normally hyperbolic
for the PRC3BP. Applying Theorem \ref{th:pert-norm-hyp} and Remark
\ref{rem:pert-norm-hyp} we obtain a family of normally hyperbolic manifolds
$\Lambda_{e}$ locally invariant under $\Phi_{t}^{e}$, and a function
$F:\Lambda\times[0,e_{0}(\mu)]\rightarrow R^{4}\times S^{1}$ such that
\begin{equation}
F(\Lambda,e)=\Lambda_{e}=\{\left(  \Lambda_{t,e},t\right)  |\Lambda
_{t,e}\subset\mathbb{R}^{4},t\in S^{1}\}.\nonumber
\end{equation}
From the Implicit Function Theorem we know that the libration point $\bar
{L}_{2}^{\mu}$ continues for small values of $e$ to a $2\pi$ periodic orbit
$\bar{L}_{2}^{\mu,e}(t)$. We can modify $F$ so that $F(\bar{L}_{2}^{\mu
},t,e)=\bar{L}_{2}^{\mu,e}(t).$ By Remark \ref{rem:pert-norm-hyp} the function
$F$ is $C^{r_{1}},$ where $r_{1}=\alpha/\beta.$ From the proof of Lemma
\ref{lem:Lap-orb normally hyp.} we know that for sufficiently small $\mu$ we
have $\alpha\thickapprox\sqrt{1+2\sqrt{7}}$ and that $\beta>0$ can be chosen
to be arbitrarily close to zero. This means that for sufficiently small $\mu,$
the function $F$ is $C^{k}$ for any given $k>0.$ Since $\Lambda_{e}$ is
locally invariant under $\Phi_{t}^{e}$ and the flow is $2\pi$ periodic, for
any $t_{0}\in S^{1}$ the manifold $\Lambda_{t_{0},e}$ is locally invariant
under $P_{t_{0}}^{e}.$

Let us fix $t_{0}=0,$ fix small $\mu>0$ and fix a Poincar\'{e} map
$P^{e}:=P_{t_{0}=0}^{e}$ (here we could consider a map $P_{t_{0}}^{e}$ for any
$t_{0}\in S^{1},$ but we fix $t_{0}=0$ for simplicity) and consider $e$ such
that $e\mu^{-2/3}<\kappa$, which ensures that (\ref{eq:He-reduced}) is valid.
We shall use a notation $B_{R}=\pi_{\{t=0\}}\{l(r)|r\leq R\}$ for the set of
Lyapounov orbits with radius smaller or equal to $R$. We will now show that
for sufficiently small $e\mu^{-1/3}$ the Poincar\'{e} map
\begin{equation}
P^{e}:F(B_{R},0,e)\rightarrow\Lambda_{0,e} \label{eq:poinc-loc}%
\end{equation}
is properly defined and symplectic. By (\ref{eq:He-reduced}), for sufficiently
small $e\mu^{-1/3}$ we can see that (\ref{eq:poinc-loc}) is properly defined
because $\Lambda_{0,e}$ is locally invariant. The map $P^{e}$ is a restriction
to $\Lambda_{0,e}$ of a time $2\pi$ shift along a trajectory of a Hamiltonian
system. Such shift is symplectic for the standard form $\omega=d\bar{x}\wedge
d\bar{p}_{x}+d\bar{y}\wedge d\bar{p}_{y}$ (here we use notation $\mathbf{\bar
{x}}=\left(  \bar{x},\bar{y},\bar{p}_{x},\bar{p}_{y}\right)  $ for coordinates
since we are working in Hill's coordinates (\ref{eq:Hill-coordinates})). In
order to show that $P^{e}$ is symplectic it is therefore sufficient to show
that $\omega$ is non degenerate on $\Lambda_{0,e}.$ For sufficiently small
$\mu$ and $e$ the manifold $\Lambda_{0,e}$ is arbitrarily close to the
manifold of the Lyapounov orbits of the Hill's problem
(\ref{eq: Hills problem}), which in turn, for $r$ sufficiently close to zero,
is arbitrarily close to the vector space $V$ given by the eigenvectors of the
pure complex eigenvalues $\pm\alpha_{2}^{\text{Hill}}=\pm\sqrt{1-2\sqrt{7}}.$
To show that $\omega$ is not degenerate on $\Lambda_{0,e}$ it is therefore
sufficient to show that $\omega$ is not degenerate on $V.$ The eigenvectors
corresponding to $\pm\alpha_{2}^{Hill}$ are $v$ and $-i\bar{v}$, where $v$ is
the second column in $\Phi^{Hill}$ (see (\ref{eq: Phi=})), which was
symplectic, therefore $\omega(v,-i\bar{v})=1$. The space $V$ is spanned by
$x_{1}$ and $x_{2}$, where $v=x_{1}+ix_{2}$. An easy computation shows that
$1=\omega(v,-i\bar{v})=-2\omega(x_{1},x_{2})$, which means that $\omega$ is
not degenerate on $V$.

Now we will use the KAM Theorem \ref{th:KAM} to show that most of the
Lyapounov orbits on $\Lambda_{0,e}$ survive under a sufficiently small
perturbation. Let us first note that even though the Theorem is stated for a
map $f$ on $[0,1]\times\mathbb{T}^{1}$, the KAM result is local by its nature
and also holds for a map $f:[0,1]\times\mathbb{T}^{1}\rightarrow
\mathbb{R}\times\mathbb{T}^{1}$, as will be the case in our setting. Let
$\omega$ denote the standard symplectic form in $\mathbb{R}^{4}$ i.e.
$\omega=d\bar{x}\wedge d\bar{p}_{x}+d\bar{y}\wedge d\bar{p}_{y}.$ Let
$\omega^{e}$ denote the induced form on $\Lambda_{0,e}.$ There exist
$C^{r_{1}-2}$ (jointly with the parameter $e$) close to identity coordinate
maps $c_{e}:\Lambda_{0,e}\rightarrow\Lambda_{0,0}$ which transport the
symplectic forms $\omega^{e}$ into the standard one (see \cite[page
367]{Llave}). The map
\begin{equation}
\bar{P}^{e}:=c_{e}\circ P^{e}\circ\left(  c_{e}\right)  ^{-1}:B_{R}\rightarrow
B_{R_{\mathrm{Hill}}}\nonumber
\end{equation}
is properly defined for sufficiently small $e$. Clearly for $e=0$ we have
$\bar{P}^{0}=P^{0}.$ From the fact that $P^{e}$ is symplectic and the fact
that $P^{0}=P_{t_{0}=0}^{e=0}$ is a twist map follow the same properties for
our maps $\bar{P}^{e}$ and $\bar{P}^{0}$ respectively. Now we pass to the
action angle coordinates. From Lemma \ref{lem: if a_2 ne 0 then P is a twist}
we have that $\bar{P}^{0}$ has the form (\ref{eq:KAM-form}). Exact
simplecticity of $\bar{P}^{e}$ is a direct consequence of simplecticity
combined with invariance of the origin. To apply the KAM Theorem \ref{th:KAM}
what is now left is to show that
\begin{equation}
\left\Vert \bar{P}^{e}-\bar{P}^{0}\right\Vert _{C^{r_{1}-2}}=O(e\mu^{-1/3}).
\label{eq:Pe-close-to-Po}%
\end{equation}
This comes from the fact that in the neighborhood of $\bar{L}_{2}^{\mu}$ the
perturbing term in (\ref{eq:He-reduced}) is uniformly $O(e\mu^{-1/3})$ in the
$C^{l}$ norm. Observe that this estimate holds both in the original
coordinates and in the action angle coordinates, because the origin is the
fixed point for $\bar{P}^{e}$. Therefore the time $2\pi$ shift maps $P_{t_{0}%
}^{e}$ and $P_{t_{0}}^{0}$ are also $O(e\mu^{-1/3})$ close, from which
(\ref{eq:Pe-close-to-Po}) follows. This gives us a Cantor family of invariant
tori $l_{0}^{e}\left(  r\right)  $ for $r\in\mathfrak{C.}$ Now for
$r\in\mathfrak{C}$ we can define
\begin{align*}
l^{e}(r)  &  :=\left\{  \Phi_{t}^{e}(x,0)|x\in l_{0}^{e}(r),t\in\lbrack
0,2\pi)\right\} \\
l_{t_{0}}^{e}(r)  &  :=\Phi_{t_{0}}^{e}(l_{0}^{e}\left(  r\right)  ,0).
\end{align*}

The fact that the complement of the Cantor set $\mathfrak{C}$ is $O(\left(
e\mu^{-1/3}\right)  ^{1/2})$ follows from the KAM Theorem (see also
\cite{Poschel} for more details).

In the above argument we require that $e\mu^{-1/3}<c$ with some
sufficiently small $c$ (which is independent of $\mu$) so that
both the normally hyperbolic theorem and KAM can be applied. Let
finish by observing that for any given $c$, by choosing
sufficiently small $\mu^{\ast}>0$ and requiring that
$\mu<\mu^{\ast}$ and $e\mu^{-2/3}<\kappa$ the estimate
$e\mu^{-1/3}<c$ follows.
\end{proof}

\begin{remark}
\label{rem:KAM-for-all-k}An identical argument to the above proof can be
performed to obtain a mirror result to Theorem \ref{th:KAM for L2} for any
fixed parameter $\mu_{k}$. To do so though one would have to verify the twist
condition. This has been verified numerically for a sequence of parameters in
Table 1. It is visible that as the masses $\mu_{k}$ decrease the twist
coefficient increases to infinity. It is therefore reasonable to believe that
the twist condition holds for all parameters $\mu_{k}$. Let us emphasize that
for sufficiently large $k$ we have rigorously verified the twist condition in
Lemma \ref{lem:twist-c}, hence the claim of Theorem \ref{th:KAM for L2} holds
for $\mu_{k}$ sufficiently large $k.$
\end{remark}

\begin{remark}
We can set $R$ from Theorem \ref{th:KAM for L2} to be the radius of one of the
surviving tori. This will ensure that $\Lambda_{e}$ is an invariant manifold
with a boundary $l^{e}(R)$. We should also emphasize that the choice of $R$ is
independent of $\mu$.
\end{remark}

If for two surviving tori $l^{e}(r_{1})$ and $l^{e}(r_{2})$ with $r_{1}<r_{2}$
there does not exist an $r\in(r_{1},r_{2})$ for which the torus $l(r)$ is
perturbed to an invariant torus of the PRE3BP, then we say that there exists a
gap between the tori $l^{e}(r_{1})$ and $l^{e}(r_{2})$.

\begin{proposition}
\label{prop:KAM-gaps}Let $e\mu^{-1/3}$ be sufficiently small so that the claim
of Theorem \ref{th:KAM for L2} holds. Then there exists an interval
$I\subset\lbrack0,R]$ with measure of order $\left(  e\mu^{-1/3}\right)
^{1/2}$, for which the set $I\cap\mathfrak{C}$ for which Lyapounov orbits
persist under perturbation has gaps smaller than $\zeta e\mu^{-1/3}$, where
$\zeta>0$ is any given constant.
\end{proposition}

\begin{proof}
The fact that such an interval exists will follow from the fact that the
complement of the Cantor set $\mathfrak{C}$ is of the measure $O(\left(
e\mu^{-1/3}\right)  ^{1/2})$. Let us divide the interval $[0,R]$ into $n$
equal parts. If on every interval the set $\mathfrak{C}$ contains gaps larger
than $\zeta e\mu^{-1/3}$, then from the fact that the measure of the
complement of $\mathfrak{C}$ is $O(\left(  e\mu^{-1/3}\right)  ^{1/2})$ (let
us say that this $O(\left(  e\mu^{-1/3}\right)  ^{1/2})$ is equal to $M\left(
e\mu^{-1/3}\right)  ^{1/2}$ for some $M>0)$ the number of such intervals $n$
must satisfy
\begin{equation}
n\zeta e\mu^{-1/3}\leq M\left(  e\mu^{-1/3}\right)  ^{1/2},\nonumber
\end{equation}
which means that $n\leq\frac{1}{\zeta}M\left(  e\mu^{-1/3}\right)  ^{-1/2}.$
If we divide the interval $[0,R]$ into a slightly larger number $\tilde{n}$ of
equal intervals then at least one of them (this will be our interval $I$)
cannot contain a gap larger than $\zeta e\mu^{-1/3}$. The size of such an
interval is equal to
\begin{equation}
\frac{R}{\tilde{n}}\approx\frac{1}{M}R\zeta\left(  e\mu^{-1/3}\right)
^{1/2}.\nonumber
\end{equation}

\end{proof}

\section{Melnikov method}

\label{sec:melnikov}In the previous section we have shown that the normally
hyperbolic manifold with a boundary $\Lambda=\left\{  l(r):r\in\lbrack
0,R]\right\}  $ (considered in the extended phase space) of the PRC3BP
(\ref{eq:He-reduced}) in Hill's coordinates (\ref{eq:PRC3BP-Hill}) persists
under perturbation to $\Lambda_{e}$ which is a normally hyperbolic invariant
manifold with a boundary of the PRE3BP with eccentricity $e.$ Moreover, we
have shown that $\Lambda_{e}$ contains a Cantor set of two dimensional
invariant KAM tori. In this section we will consider the problem of
intersections of the stable and unstable manifolds of such tori.

In this section we shall once again work in the Hill's coordinates
(\ref{eq:Hill-coordinates}) and Hamiltonian (\ref{eq:PRE3BP-Hill}). It will
also be convenient for us to parameterize the manifolds $\Lambda$ and
$\Lambda_{e}$ using the radius angle coordinates $r,\varphi$ of the Lyapounov
orbits from Section \ref{sec:lap-twist} together with time $t\in S^{1}$. For
$e=0$ we thus use a natural parameterization of the Lyapounov orbits by their
Birkhoff normal form coordinates (coordinates obtained from Theorem
\ref{th: Twierdzenie Mosera}). After the perturbation it will be enough for us
to use the fact that the parametrization is smooth (in fact, from the proof of
Theorem \ref{th:KAM for L2} we know that it will be $C^{r_{1}}$ with
$r_{1}=\alpha/\beta$) and that we can parametrize (see proof of Theorem
\ref{th:KAM for L2}) the perturbed libration point $\bar{L}_{2}^{\mu}%
=\mu^{-1/3}\left(  L_{2}^{\mu}-\left(  \mu-1,0,0,\mu-1\right)  \right)  $
(which continues to a $2\pi$ periodic orbit) by $r=0.$ We can not assume
though that the surviving perturbed KAM tori are parameterized by $r$. Our
parametrization simply follows from the Normally Hyperbolic Invariant Manifold
Theorem (see Theorem \ref{th:pert-norm-hyp} and Remarks
\ref{rem:pert-norm-hyp}, \ref{rem:pert-norm-hyp-with-manifolds}) without
involving the KAM Theorem.

Let us recall that for $\mu$ close to $\mu_{k}$ from Theorem
\ref{lem homoclinic orbit for uk} prior to the perturbation the fibres
$W_{p}^{s}$ for $p\in\Lambda$ intersect transversally with the section
$\{\bar{y}=0\}$ (see Section \ref{sec:prelim}, Theorems
\ref{lem homoclinic orbit for uk}, \ref{lem: transversality-Simo} and also
\cite{Simo}). The same goes for the unstable fibers $W_{p}^{u}$. This means
that for sufficiently small $0<e$ the same will hold for the stable fibres
$W_{p}^{s,e}$ and unstable fibres $W_{p}^{u,e}$ of points $p\in\Lambda_{e}.$
Each point $p\in\Lambda_{e}$ can be parameterized by $\mu,r,\varphi$ and
$t_{0}.$ The fibres of such points $W_{p}^{u,e},$ $W_{p}^{s,e}$ are one
dimensional and contained in sections $\Sigma_{t_{0}}=\{(q,t_{0}%
)|q\in\mathbb{R}^{4}\}.$ The intersections of $W_{p}^{u,e}$ and of
$W_{p}^{s,e}$ with $\{\bar{y}=0\}$ are functions of $(\mu,r,\varphi,t_{0},e)$.
For a point $p\in\Lambda_{e}$ parameterized by $\mu,r,\varphi$ and $t_{0}$ we
introduce the following notation for the first intersections of $W_{p}^{s,e}$
and of $W_{p}^{u,e}\ $with $\{\bar{y}=0\}$
\begin{align*}
(p^{s}(r,\varphi,t_{0},e),t_{0})  &  =W_{p}^{s,e}\cap\{y=0\},\\
(p^{u}(r,\varphi,t_{0},e),t_{0})  &  =W_{p}^{u,e}\cap\{y=0\}.
\end{align*}
Let $q^{s}(r,\varphi,t_{0},e,t)$ and $q^{u}(r,\varphi,t_{0},e,t)$ be the
orbits (considered in the standard (not extended) phase space) of the PRE3BP,
which start from the points $p^{s}(r,\varphi,t_{0},e)$ and $p^{u}%
(r,\varphi,t_{0},e)$ respectively at time $t=t_{0}$ i.e.%
\begin{align*}
q^{s}(r,\varphi,t_{0},e,t_{0})  &  =p^{s}(r,\varphi,t_{0},e),\\
q^{u}(r,\varphi,t_{0},e,t_{0})  &  =p^{u}(r,\varphi,t_{0},e).
\end{align*}
($p^{s},p^{u},q^{s},q^{u}$ depend on the choice of $\mu$, but we omit this in
our notations for simplicity). Let us note that for parameters $\mu=\mu_{k}$
from Theorem \ref{lem homoclinic orbit for uk}
\begin{align}
q^{s}(0,\varphi,t_{0},0,t)  &  =q^{u}(0,\varphi,t_{0},0,t)\label{eq:qs=q0}\\
&  =\bar{q}_{\mu_{k}}^{0}(t-t_{0})\nonumber\\
&  =\mu^{-1/3}\left(  q_{\mu_{k}}^{0}(t-t_{0})-\left(  \mu_{k}-1,0,0,\mu
_{k}-1\right)  \right)  ,\nonumber
\end{align}
where $q_{\mu_{k}}^{0}(t)$ is the homoclinic orbit to $L_{2}^{\mu_{k}}$ in the
PRC3BP defined just after the statement of
Theorem~\ref{lem homoclinic orbit for uk}. Let us remind the reader that
$q_{\mu_{k}}^{0}(0)\in\{y=0\}$.

\begin{lemma}
\label{lem:variational}For fixed $\mu$ and $i\in\{s,u\}$
\begin{align}
q^{i}(r,\varphi,t_{0},0,t_{0})  &  =q^{i}(0,\varphi,t_{0},0,t_{0}%
)+O(r),\label{eq:q-exp1}\\
q^{i}(r,\varphi,t_{0},e,t_{0})  &  =q^{i}(r,\varphi,t_{0},0,t_{0}%
)+e\frac{\partial q^{i}}{\partial e}(r,\varphi,t_{0},0,t_{0}%
)+o(e),\label{eq:q-expansion}\\
\frac{\partial q^{i}}{\partial e}(r,\varphi,t_{0},0,t_{0})  &  =\frac{\partial
q^{i}}{\partial e}(0,\varphi,t_{0},0,t_{0})+O(r), \label{eq:q-exp3}%
\end{align}
where the bounds $o(e)$ and $O(r)$ are independent from $t_{0}$.

In addition for $f,g$ from (\ref{eq:PRE3BP-reduced})
\begin{align}
\frac{d}{dt}\left(  \frac{\partial q^{i}}{\partial e}(0,\varphi,t_{0}%
,0,t)\right)   &  =Df\left(  \mu,q^{i}(0,\varphi,t_{0},0,t)\right)
\frac{\partial q^{i}}{\partial e}(0,\varphi,t_{0},0,t)\label{eq:variational}\\
&  +g\left(  q^{i}(0,\varphi,t_{0},0,t),t\right)  ,\nonumber
\end{align}
and $\frac{\partial q^{i}}{\partial e}(0,\varphi,t_{0},0,t)$ is bounded for
all $t\in[ t_{0},+\infty)$ for $i=s$ (or $t\in(-\infty,t_{0}]$ for $i=u$).
\end{lemma}

\begin{proof}
The normally hyperbolic manifold and the foliation of its stable and unstable
manifolds behave smoothly under perturbation and equations (\ref{eq:q-exp1}%
--\ref{eq:q-exp3}) are a simple consequence of this.

It remains to prove (\ref{eq:variational}). Let $i=s$. We have
\begin{equation}
\frac{d}{dt}\left(  \frac{\partial q^{s}}{\partial e}(0,\varphi,t_{0}%
,0,t)\right)  =\frac{\partial}{\partial e}\frac{d}{dt}q^{s}(0,\varphi
,t_{0},0,t)\quad\quad\quad\quad\quad\quad\quad\nonumber
\end{equation}%
\begin{align}
&  =\frac{\partial}{\partial e}\left(  f\left(  \mu,q^{s}(0,\varphi
,t_{0},e,t)\right)  +eg\left(  q^{s}(0,\varphi,t_{0},e,t),t\right)  +O(\left(
e\mu^{-1/3}\right)  ^{2})\right)  |_{e=0}\nonumber\\
&  =Df\left(  \mu,q^{s}(0,\varphi,t_{0},0,t)\right)  \frac{\partial q^{s}%
}{\partial e}(0,\varphi,t_{0},0,t)+g\left(  q^{s}(0,\varphi,t_{0},0,t)\right)
.\nonumber
\end{align}

The points $p^{s}(0,\varphi,t_{0},e)$ lie on stable fibres of the periodic
orbit perturbed from $L_{2}^{\mu}$. The points $p^{s}(0,\varphi,t_{0},e)$ and
$p^{s}(0,\varphi,t_{0},0)$ are therefore $O(e)$ close. Also the periodic orbit
perturbed from $L_{2}^{\mu}$ lies $O(e)$ close to $L_{2}^{\mu}.$ Since orbits
$q^{s}(0,\varphi,t_{0},e,t)$ start from $p^{s}(0,\varphi,t_{0},e)$ this gives
us
\begin{equation}
\left\vert q^{s}(0,\varphi,t_{0},e,t)-q^{s}(0,\varphi,t_{0},0,t)\right\vert
=O(e) \label{eq:qe-bound}%
\end{equation}
for $t\in\lbrack t_{0},+\infty).$ This means that%
\begin{equation}
\left\vert \frac{\partial q^{s}}{\partial e}(0,\varphi,t_{0},0,t)\right\vert
=\left\vert \lim_{e\rightarrow0}\frac{q^{s}(0,\varphi,t_{0},e,t)-q^{s}%
(0,\varphi,t_{0},0,t)}{e}\right\vert \nonumber
\end{equation}
is bounded.

For $i=u$ the argument is analogous.
\end{proof}

\begin{remark}
\label{rem:radius-partial-info}Let us note that the bound $O(r)$
in (\ref{eq:q-exp1}) is independent of $\mu$.
\end{remark}

\begin{proof}
This follows from formulas for the parameterization of the
intersection of the manifolds with $\{y=0\}$ from Theorem
\ref{lem: transversality-Simo}. Each curve from Theorem \ref{lem:
transversality-Simo} gives an intersection of an unstable manifold
of a Lyapounov orbit. A point $q^{i}(r,\varphi,t_{0},0,t_{0})$ is
a point of intersection of the unstable manifold of a Lyapounov
orbit $l(r)$ with $\{y=0\},$ and in $x,\dot{x}$ coordinates is
represented by a point $\left( x(\sigma,\sqrt{\Delta
C}),\dot{x}(\sigma,\sqrt{\Delta C})\right)  $ on a curve from
Theorem \ref{lem: transversality-Simo}, for some $\sigma\in S^{1}$
and $\Delta C>0.$ The point $q^{i}(0,\varphi,t_{0},0,t_{0})$ is
the point of intersection of the homoclinic orbit to
$L_{2}^{\mu_{k}}$ with $\{y=0\}$ and
in $x,\dot{x}$ coordinates is given by $\left(  x(\sigma,0),\dot{x}%
(\sigma,0)\right)  $ (where the choice of $\sigma$ plays no role
since for $\Delta C=0$ equations from Theorem \ref{lem:
transversality-Simo} give a single point).

From Theorem \ref{lem: transversality-Simo}
\begin{equation}
x(\sigma,\sqrt{\Delta C})-x(\sigma,0)=O(\sqrt{\Delta C})\quad\text{and\quad
}\dot{x}(\sigma,\sqrt{\Delta C})-\dot{x}(\sigma,0)=O(\sqrt{\Delta
C}).\label{eq:x-dC-distance}%
\end{equation}
By (\ref{eq:PRC3BP-Hill}) and Lemma \ref{lem:dist(Lc,L2)}
\begin{align}
\sqrt{\Delta C}  &  =\sqrt{H(\mu,\mu^{1/3}l(r)+\left(  \mu-1,0,0,\mu-1\right)
)-H(\mu,L_{2}^{\mu})}\nonumber\\
&  =\sqrt{\mu^{2/3}\bar{H}(\mu,l(r))-\mu^{2/3}\bar{H}(\mu,l(0))}\nonumber\\
&  =\mu^{1/3}O(r).\label{eq:dC-r-relation}%
\end{align}
The points $x(\sigma,\sqrt{\Delta C}),$ $x(\sigma,0),$ $\dot{x}(\sigma
,\sqrt{\Delta C}),$ $\dot{x}(\sigma,0)$ are considered in original coordinates
of the system on the section $\{y=0\}.$ For a point $(x,\dot{x})=\left(
x(\sigma,\sqrt{\Delta C}),\dot{x}(\sigma,\sqrt{\Delta C})\right)  $ on the
section $\{y=0\}$ by (\ref{eq:H-PRC3BP}) we have $p_{x}=\dot{x}-y=\dot{x}$
and
\begin{equation}
p_{y}=\sqrt{2\left(  H(L_{2}^{\mu_{k}})+\Delta C+\Omega(x,0)\right)
-p_{x}^{2}}=\sqrt{2H(L_{2}^{\mu_{k}})}+\mu^{1/3}O(r).\nonumber
\end{equation}

Recall that the points $q^{i}(r,\varphi,t_{0},0,t_{0}),$
$q^{i}(0,\varphi,t_{0},0,t_{0})$ are given in Hill's coordinates
(\ref{eq:Hill-coordinates}). By (\ref{eq:Hill-coordinates}),
(\ref{eq:x-dC-distance}), (\ref{eq:dC-r-relation})
\begin{align*}
&  \left\vert q^{i}(r,\varphi,t_{0},0,t_{0})-q^{i}(0,\varphi,t_{0}%
,0,t_{0})\right\vert \\
&  =\mu^{-1/3}\left\vert \left(  x,0,p_{x},p_{y}\right)  (\sigma,\sqrt{\Delta
C})-\left(  x,0,p_{x},p_{y}\right)  (\sigma,0)\right\vert =O(r).
\end{align*}

\end{proof}

\begin{remark}
We believe that with techniques similar to the ones used for the
proof of Theorem \ref{lem: transversality-Simo} in \cite{Simo},
but involving additionally terms coming from the perturbation from
the PRC3BP to the PRE3BP, it should be possible to prove that the
terms $o(e)$ and $O(r)$ from (\ref{eq:q-expansion}),
(\ref{eq:q-exp3}) can be chosen independently of $\mu_{k}$. Such
statement requires a detailed proof which is rather technical. We
skip this intentionally since in later parts of our argument it
will turn out that even if this had been done by us, we still
cannot obtain uniform bounds for the radius of set on which we
have structural stability for the PRE3BP. This is due to the fact
that we have not obtained uniform bounds for the Melnikov integral
(\ref{eq:Meln-int-}) (see Remark \ref{rem:size-of-R}). Such bounds
seem even harder to obtain than proving that the terms $o(e)$ and
$O(r)$ from (\ref{eq:q-expansion}), (\ref{eq:q-exp3}) are
independent of $\mu_{k}$.

Remark \ref{rem:radius-partial-info} will allow us though to
extract some information as to the radius for which we shall have
structural stability for the PRE3BP. It will turn out that this
radius needs to converge to zero
with $\mu_{k}$ going to zero at least fast enough so that $\mu_{k}^{-1/3}%
R(\mu_{k})$ is bounded. From our proof though it is not
transparent how small exactly it shall need to be chosen.
\end{remark}

We now have the following lemma regarding the energy of the points
$p^{s}(r,\varphi,t_{0},e)$ and $p^{u}(r,\varphi,t_{0},e).$

\begin{lemma}
\label{lem:Melnikov}Assume that $\mu$ is one of the parameters $\mu_{k}$ for
which a homoclinic orbit $\bar{q}_{\mu_{k}}^{0}(t)=\mu^{-1/3}\left(
q_{\mu_{k}}^{0}(t)-\left(  \mu_{k}-1,0,0,\mu_{k}-1\right)  \right)  $ to
$\bar{L}_{2}^{\mu}=\mu^{-1/3}\left(  L_{2}^{\mu}-\left(  \mu_{k}-1,0,0,\mu
_{k}-1\right)  \right)  $ exists (See Theorem
\ref{lem homoclinic orbit for uk}). For any two points $p_{1},p_{2}\in
\Lambda_{e}$ with coordinates $(r_{1},\varphi_{1},t_{0})$ and $(r_{2}%
,\varphi_{2},t_{0})$ respectively, we have
\begin{multline*}
\bar{H}(\mu_{k},p^{s}(r_{1},\varphi_{1},t_{0},e))-\bar{H}(\mu_{k},p^{u}%
(r_{2},\varphi_{2},t_{0},e))\\
=\bar{H}(\mu_{k},l(r_{1}))-\bar{H}(\mu_{k},l(r_{2}))+eM_{\mu_{k}}(t_{0})\\
+O(e\max\left\{  |r_{1}|,|r_{2}|\right\}  )+o(e),
\end{multline*}
where%
\begin{equation}
M_{\mu_{k}}(t_{0})=\int_{-\infty}^{+\infty}\{\bar{H},\bar{G}\}(\mu_{k},\bar
{q}_{\mu_{k}}^{0}(t-t_{0}),t)dt. \label{eq:Meln-int-}%
\end{equation}

\end{lemma}

\begin{proof}
Let $\cdot$ denote the scalar product and let $\Delta_{s}$ and $\Delta_{u}$
denote the following functions
\begin{align*}
\Delta_{s}(t,t_{0})  &  :=\nabla\bar{H}(\mu_{k},\bar{q}_{\mu_{k}}^{0}%
(t-t_{0}))\cdot\frac{\partial q^{s}}{\partial e}(0,\varphi_{1},t_{0},0,t)\\
\Delta_{u}(t,t_{0})  &  :=\nabla\bar{H}(\mu_{k},\bar{q}_{\mu_{k}}^{0}%
(t-t_{0}))\cdot\frac{\partial q^{u}}{\partial e}(0,\varphi_{2},t_{0},0,t).
\end{align*}
Using the facts that $\bar{H}=\bar{H}_{r}=\bar{H}(\mu_{k},l(r))$ is constant
along the solutions $q^{s}(r,\varphi,t_{0},0,t)$ of the PRC3BP, from
(\ref{eq:q-expansion}) and (\ref{eq:qs=q0}) we can compute%
\begin{align}
\bar{H}(\mu_{k},q^{s}  &  (r_{1},\varphi_{1},t_{0},e,t_{0}))\nonumber\\
=  &  \bar{H}(\mu_{k},q^{s}(r_{1},\varphi_{1},t_{0},0,t_{0}))\nonumber\\
&  +e\nabla\bar{H}(\mu_{k},q^{s}(r_{1},\varphi_{1},t_{0},0,t_{0}))\cdot
\frac{\partial q^{s}}{\partial e}(r_{1},\varphi_{1},t_{0},0,t_{0}%
)+o(e)\nonumber\\
=  &  \bar{H}_{r_{1}}+e\nabla\bar{H}(\mu_{k},\bar{q}_{\mu_{k}}^{0}%
(0)+O(r_{1}))\cdot\left(  \frac{\partial q^{s}}{\partial e}(0,\varphi
_{1},t_{0},0,t_{0})+O(r_{1})\right)  +o(e)\nonumber\\
=  &  \bar{H}_{r_{1}}+e\nabla\bar{H}(\mu_{k},\bar{q}_{\mu_{k}}^{0}%
(0))\cdot\frac{\partial q^{s}}{\partial e}(0,\varphi_{1},t_{0},0,t_{0}%
)+O(er_{1})+o(e)\nonumber\\
=  &  \bar{H}_{r_{1}}+e\Delta_{s}(t_{0},t_{0})+O(er_{1})+o(e),
\label{eq:H-exp-1}%
\end{align}
and similarly one can show that
\begin{equation}
\bar{H}(\mu_{k},q^{u}(r_{2},\varphi_{2},t_{0},e,t_{0}))=\bar{H}_{r_{2}%
}+e\Delta_{u}(t_{0},t_{0})+O(er_{2})+o(e). \label{eq:H-exp-2}%
\end{equation}
Let us stress that, by Lemma \ref{lem:variational} we know that $O(er_{1}),$
$O(er_{2})$ and $o(e)$ are uniform with respect to $t_{0}$.

Let us investigate the evolution of $\Delta_{s}(t,t_{0})$ and $\Delta
_{u}(t,t_{0})$ in time. Let us concentrate on the term $\Delta_{s}(t,t_{0}).$
Using (\ref{eq:variational}), (\ref{eq:qs=q0}) and $\nabla\bar{H}=-Jf$ (see
(\ref{eq:def-of-f})) we can compute
\begin{align}
-\frac{d}{dt}\left(  \Delta_{s}(t,t_{0})\right)   &  =\left(  JDf(\mu_{k}%
,\bar{q}_{\mu_{k}}^{0}(t-t_{0}))\frac{d}{dt}\bar{q}_{\mu_{k}}^{0}%
(t-t_{0})\right)  \cdot\frac{\partial q^{s}}{\partial e}(0,\varphi_{1}%
,t_{0},0,t)\nonumber\\
&  \quad+\left(  Jf(\mu_{k},\bar{q}_{\mu_{k}}^{0}(t-t_{0}))\right)  \cdot
\frac{d}{dt}\frac{\partial q^{s}}{\partial e}(0,\varphi_{1},t_{0}%
,0,t)\label{eq:Melnikov-s}\\
&  =\left(  JDf(\mu_{k},\bar{q}_{\mu_{k}}^{0}(t-t_{0}))f(\mu_{k},\bar{q}%
_{\mu_{k}}^{0}(t-t_{0}))\right)  \cdot\frac{\partial q^{s}}{\partial
e}(0,\varphi_{1},t_{0},0,t)\nonumber\\
&  \quad+\left(  Jf(\mu_{k},\bar{q}_{\mu_{k}}^{0}(t-t_{0}))\right)
\cdot\left(  Df\left(  \mu_{k},\bar{q}_{\mu_{k}}^{0}(t-t_{0})\right)
\frac{\partial q^{s}}{\partial e}(0,\varphi_{1},t_{0},0,t)\right) \nonumber\\
&  \quad+\left(  Jf(\mu_{k},\bar{q}_{\mu_{k}}^{0}(t-t_{0}))\right)  \cdot
g(\mu_{k},\bar{q}_{\mu_{k}}^{0}(t-t_{0}),t)\nonumber\\
&  =\left(  Jf(\mu_{k},\bar{q}_{\mu_{k}}^{0}(t-t_{0}))\right)  \cdot g(\mu
_{k},\bar{q}_{\mu_{k}}^{0}(t-t_{0}),t)\nonumber\\
&  =-\{\bar{H},\bar{G}\}(\mu_{k},\bar{q}_{\mu_{k}}^{0}(t-t_{0}),t),\nonumber
\end{align}
where the third equality comes from the fact that for any $p,q\in
\mathbb{R}^{4}$%
\begin{equation}
\left(  JDf(\mu_{k},\bar{q}_{\mu_{k}}^{0}(t-t_{0}))p\right)  \cdot
q+(Jp)\cdot\left(  Df\left(  \mu_{k},\bar{q}_{\mu_{k}}^{0}(t-t_{0})\right)
q\right)  =0, \label{eq: JDfp*q+Jp*Dfq=0}%
\end{equation}
with $p=f(\mu_{k},\bar{q}_{\mu_{k}}^{0}(t-t_{0}))$ and $q=\frac{\partial
q^{s}}{\partial e}(0,\varphi_{1},t_{0},0,t).$ Equation
(\ref{eq: JDfp*q+Jp*Dfq=0}) follows from the fact that $\omega(p,q)=Jp\cdot q$
is the standard symplectic form which is invariant under the flow $\phi(t,x)$
of (\ref{eq:PRC3BP-Hill}) i.e.
\begin{equation}
\omega\left(  \frac{\partial}{\partial x}\phi(t,x)p,\frac{\partial}{\partial
x}\phi(t,x)q\right)  =\omega(p,q), \label{eq:form-invariant}%
\end{equation}
hence by differentiating (\ref{eq:form-invariant}) with respect to $t$ and
setting $t=0$ we obtain (\ref{eq: JDfp*q+Jp*Dfq=0}).

We can now compute $\Delta_{s}(t_{0},t_{0})$ using (\ref{eq:Melnikov-s})%
\begin{equation}
\Delta_{s}(+\infty,t_{0})-\Delta_{s}(t_{0},t_{0})=\int_{t_{0}}^{+\infty}%
\{\bar{H},\bar{G}\}(\mu_{k},\bar{q}_{\mu_{k}}^{0}(t-t_{0}),t)dt.\nonumber
\end{equation}
Since $\lim_{t\rightarrow+\infty}\bar{q}_{\mu_{k}}^{0}(t-t_{0})=\bar{L}%
_{2}^{\mu_{k}}$ at geometric rate and $f(\mu_{k},\bar{L}_{2}^{\mu_{k}})=0$,
from the fact that $\frac{\partial q^{s}}{\partial e}(0,\varphi_{1}%
,t_{0},0,t)$ is bounded on $[t_{0},+\infty)$ we have
\begin{equation}
\Delta_{s}(+\infty,t_{0})=\lim_{t\rightarrow+\infty}Jf(\mu_{k},\bar{q}%
_{\mu_{k}}^{0}(t-t_{0}))\cdot\frac{\partial q^{s}}{\partial e}(0,\varphi
_{1},t_{0},0,t)=0\nonumber
\end{equation}
and therefore
\begin{equation}
-\Delta_{s}(t_{0},t_{0})=\int_{t_{0}}^{+\infty}\{\bar{H},\bar{G}\}(\mu
_{k},\bar{q}_{\mu_{k}}^{0}(t-t_{0}),t)dt, \label{eq:int-1}%
\end{equation}
and $\Delta$ is uniformly with respect to $t_{0}$ absolutely convergent.

Analogous computations give
\begin{equation}
\Delta_{u}(t_{0},t_{0})=\int_{-\infty}^{t_{0}}\{\bar{H},\bar{G}\}(\bar{q}%
_{\mu_{k}}^{0}(t-t_{0}),t)dt. \label{eq:int-2}%
\end{equation}
From (\ref{eq:H-exp-1}), (\ref{eq:H-exp-2}), (\ref{eq:int-1}) and
(\ref{eq:int-2}) we obtain our claim.
\end{proof}

\begin{theorem}
\label{th:Melnikov} Consider the PRE3BP with a sufficiently small parameter
$\mu=\mu_{k}$ for which a homoclinic orbit to $L_{2}^{\mu_{k}}$ exists (See
Theorem \ref{lem homoclinic orbit for uk} ). Assume that
\begin{equation}
M_{\mu_{k}}(t_{0})=\int_{-\infty}^{+\infty}\{\bar{H},\bar{G}\}(\mu_{k},\bar
{q}_{\mu_{k}}^{0}(t-t_{0}),t)dt \label{eq:Melnikov}%
\end{equation}
has simple zeros. Let $R(\mu_{k})\in\mathbb{R}$, be such that $0<R(\mu_{k})$
and $\mu_{k}^{-1/3}R(\mu_{k})$ is sufficiently close to zero. Then for any
$R\in(0,R(\mu_{k}))$ there exists an $e_{0}(R)$ such that for all $e\mu
_{k}^{-1/3}\in[ 0,\min(e_{0}(R),\kappa\mu_{k}^{1/3})]$ (where $\kappa$ is the
constant from (\ref{eq:kappa-e-bound} )) and all $r\in\mathfrak{C}\cap[
R,R(\mu_{k})]$ for which $l(r)$ is perturbed to an invariant torus
$l^{e}\left(  t\right)  $
\begin{equation}
l^{e}(r)=\{(l_{t_{0}}^{e}(r),t_{0})|t_{0}\in S^{1}\}\nonumber
\end{equation}
the manifolds $W_{l^{e}(r)}^{s}$ and $W_{l^{e}(r)}^{u}$ (considered in the
extended phase space) intersect transversally.
\end{theorem}

\begin{proof}
We consider the PRE3BP (\ref{eq:PRE3BP-reduced}) in Hill's coordinates
(\ref{eq:Hill-coordinates}). In \cite{Simo} it has been shown (see also
Theorem \ref{lem: transversality-Simo} in Section \ref{subsec:PRC3BP}) that
for the unperturbed PRC3BP $W_{l(r)}^{s}$ intersects transversally with
$W_{l(r)}^{u}$ at $\{\bar{y}=0\}.$ Let $\Sigma_{t_{0}}=\mathbb{R}^{4}%
\times\{t_{0}\}$ be the time $t_{0}$ section in the extended phase space. Let
$v^{0}(t_{0})$ denote some point for which%
\begin{equation}
v^{0}(t_{0})\in W_{l(r)}^{s}\cap W_{l(r)}^{u}\cap\{\bar{y}=0\}\cap
\Sigma_{t_{0}}. \label{eq:voto}%
\end{equation}
There can be more than just one such point (see Figures \ref{fig:tubes},
\ref{fig:p(c,e)}), namely, if we consider $\pi_{x,p_{x}}(W_{l(r)}^{u}%
\cap\Sigma_{t_{0}}\cap\{\bar{y}=0\})$ and $\pi_{\bar{x},\bar{p}_{x}}%
(W_{l(r)}^{s}\cap\Sigma_{t_{0}}\cap\{y=0\})$ then the two sets are
homeomorphic to two circles, which intersect transversally at least one point
$(\bar{x}=\bar{x}_{0},\bar{p}_{x}=0)$ (see Theorem
\ref{lem: transversality-Simo}). Fixing any one of such points of intersection
will be sufficient for our proof. From the construction we have $\pi
_{x,y,p_{x},t}(v^{0}(t_{0}))=(\bar{x}_{0},0,0,t_{0}).$ The extended phase
space of the PRC3BP is five dimensional. We can choose the energy $\bar{H}$ to
be the remaining fifth coordinate in the neighborhood of $v^{0}(t_{0}).$ In
these coordinates%
\begin{equation}
v^{0}(t_{0})=\left(  \bar{x}=\bar{x}_{0},\bar{y}=0,\bar{p}_{x}=0,\bar{H}%
=\bar{H}(\mu_{k},l(r)),t=t_{0}\right)  .\nonumber
\end{equation}
\begin{figure}[ptb]
\begin{center}
\includegraphics[
height=1.3in
]{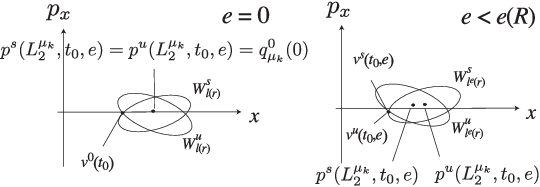}
\end{center}
\caption{The manifolds $W_{l(r)}^{s},$ $W_{l(r)}^{u},$ $W_{l^{e}(r)}^{s}$ and
$W_{l^{e}(r)}^{u}$ for $r>R,$ intersected with $\Sigma_{t_{0}}$ and $\{y=0\},$
and projected onto the $x,p_{x}$ coordinates.}%
\label{fig:p(c,e)}%
\end{figure}

Now we perturb from the PRC3BP to the PRE3BP. Consider sufficiently small
$\mu_{k}$ and perturbation $e$ satisfying $e\mu_{k}^{-2/3}<\kappa$ so that by
Theorem \ref{th:KAM for L2} we have a Cantor set $\mathfrak{C}$ of KAM tori
$l^{e}(r)$\. Let us consider some small $R(\mu_{k})\in\mathbb{R}$ (the size of
$R(\mu_{k})$ needs to be chosen small compared with $M_{\mu_{k}}(t)$. This is
discussed later on in our argument). Let $0<R<R(\mu_{k})$ and $r\in
\mathfrak{C}\cap\lbrack R,R(\mu_{k})].$ Consider now the following sets%
\begin{equation}
\pi_{\bar{x},\bar{p}_{x}}\left(  W_{l^{e}(r)}^{u}\cap\Sigma_{t_{0}}\cap
\{\bar{y}=0\}\right)  \quad\text{and\quad}\pi_{\bar{x},\bar{p}_{x}}\left(
W_{l^{e}(r)}^{s}\cap\Sigma_{t_{0}}\cap\{\bar{y}=0\}\right)  .\nonumber
\end{equation}
For sufficiently small $e$ these sets remain homeomorphic to circles. For
$e=0$ the curves intersect transversally at an angle $O(\mu_{k}^{1/3})$ (see
Remark \ref{rem:curve-splitting}). This means that by choosing $e\mu^{-1/3}>0$
sufficiently small the curves shall intersect at some point $(\bar{x}=\bar
{x}_{0}(t_{0},e),\bar{p}_{x}=\bar{p}_{x}^{0}(t_{0},e))$ which is close to
$(x_{0},0)$. The choice of $e\mu^{-1/3}$ also depends on $R$, since in order
to ensure that the curves intersect we assume that their radius is greater
than $R$. Hence we have $e_{0}(R)$ in the formulation of our theorem. For the
PRE3BP the energy $\bar{H}$ is no longer preserved. This means that the
intersection of the circles on the $\bar{x},\bar{p}_{x}$ plane does not imply
an intersection in the extended phase space. Namely we have two points
$v^{s}(t_{0},e)\in W_{l^{e}(r)}^{s}\cap\Sigma_{t_{0}}\cap\{\bar{y}=0\}$ and
$v^{u}(t_{0},e)\in W_{l^{e}(r)}^{u}\cap\Sigma_{t_{0}}\cap\{\bar{y}=0\}$ which
may differ on the energy coordinate (see Figure \ref{fig:p(c,e)-3d})%
\begin{align}
v^{s}(t_{0},e)  &  =\left(  \bar{x}_{0}(t_{0},e),0,\bar{p}_{x}^{0}%
(t_{0},e),h^{s}(t_{0},e),t_{0}\right) \label{eq:hu-hs}\\
v^{u}(t_{0},e)  &  =\left(  \bar{x}_{0}(t_{0},e),0,\bar{p}_{x}^{0}%
(t_{0},e),h^{u}(t_{0},e),t_{0}\right)  .\nonumber
\end{align}
We will show that assumptions of our theorem imply that for some $t_{0}$ the
points $v^{s}(t_{0},e)$ and $v^{u}(t_{0},e)$ coincide. This will imply
intersection between $W_{l^{e}(r)}^{s}$ and $W_{l^{e}(r)}^{u}$. Later we will
also show that such intersection is transversal.

The points $v^{s}(t_{0},e)$ and $v^{u}(t_{0},e)$ are both contained in
$\{\bar{y}=0\}.$ Moreover $v^{s}(t_{0},e)\in W_{l^{e}(r)}^{s}\cap\Sigma
_{t_{0}}=W_{l_{t_{0}}^{e}(r)}^{s}$ and $v^{u}(t_{0},e)\in W_{l^{e}(r)}^{u}%
\cap\Sigma_{t_{0}}=W_{l_{t_{0}}^{e}(r)}^{u}.$ This means that there exist
$r^{s},$ $\varphi^{s}$ and $r^{u},$ $\varphi^{u}$ (these depend on $t_{0}$ and
$e$ but we omit this in our notations for simplicity) such that%
\begin{align}
v^{s}(t_{0},e)  &  =(p^{s}(r^{s},\varphi^{s},t_{0},e),t_{0})\label{eq:ru-rs}\\
v^{u}(t_{0},e)  &  =(p^{u}(r^{u},\varphi^{u},t_{0},e),t_{0}).\nonumber
\end{align}
In the proof of Theorem \ref{th:KAM for L2} $l_{t_{0}}^{e}(r)$ is constructed
from continuation along trajectories of a KAM torus $l_{0}^{e}(r).$ Therefore
from (\ref{eq:kam-tor-close}) in KAM Theorem \ref{th:KAM} we have that%
\begin{align}
r^{s}  &  =r+O(e\mu^{-1/3}),\label{eq:r-tori-bound}\\
r^{u}  &  =r+O(e\mu^{-1/3}),\nonumber
\end{align}
and the bound $O(e\mu^{-1/3})$ is uniform for all $r\in\mathfrak{C}$. Applying
Lemma \ref{lem:Melnikov} we have that%
\begin{align}
h^{u}(t_{0},e)-h^{s}(t_{0},e)  &  =\bar{H}(\mu_{k},p^{s}(r^{s},\varphi
^{s},t_{0},e))-\bar{H}(\mu_{k},p^{u}(r^{u},\varphi^{u},t_{0}%
,e))\label{eq:Meln-tori}\\
&  =\bar{H}(\mu_{k},l(r^{s}))-\bar{H}(\mu_{k},l(r^{u}))+eM_{\mu_{k}}%
(t_{0})+O(eR(\mu_{k}))+o(e).\nonumber
\end{align}
The Hamiltonian in coordinates $(q,p)=(q_{1},q_{2},p_{1},p_{2})=(\bar{x}%
,\bar{y},\bar{p}_{x},\bar{p}_{y})-\bar{L}_{2}^{\mu_{k}}$ centered in $\bar
{L}_{2}^{\mu_{k}}$ is%
\begin{equation}
H\left(  p,q\right)  :=\bar{H}\left(  \mu_{k},(q,p)+\bar{L}_{2}^{\mu_{k}%
}\right)  .\nonumber
\end{equation}
By Lemma \ref{lem:dist(Lc,L2)} we hence know that
\begin{equation}
\bar{H}(\mu_{k},l(r))=\bar{H}(\mu_{k},\bar{L}_{2}^{\mu_{k}})+\frac{1}{2}%
D^{2}H(0)\left(  \Phi(0,1,0,i)\right)  r^{2}+o(r^{2})\nonumber
\end{equation}
hence from (\ref{eq:r-tori-bound}) and (\ref{eq:Meln-tori}) we have%
\begin{align}
h^{u}(t_{0},e)-h^{s}(t_{0},e)  &  =O(re\mu^{-1/3})+eM_{\mu_{k}}(t_{0}%
)+O(eR(\mu_{k}))+o(e)\label{eq:H-difference}\\
&  =eM_{\mu_{k}}(t_{0})+O(e\mu^{-1/3}R(\mu_{k}))+o(e).\nonumber
\end{align}
Setting first $R(\mu_{k})$ sufficiently small (so that $\mu_{k}^{-1/3}%
R(\mu_{k})$ is small in comparison to $M_{\mu_{k}} (t)$) and then reducing $e$
sufficiently close to zero implies that since $M_{\mu_{k}}(t_{0})$ has simple
zeros, for some parameters $t_{0}$ (close to these zeros) we will have
$h^{u}(t_{0},e)-h^{s}(t_{0},e)=0,$ which implies that $v^{s}(t_{0}%
,e)=v^{u}(t_{0},e)$ and in turn ensures that%
\begin{equation}
W_{l^{e}(r)}^{s}\cap W_{l^{e}(r)}^{u}\neq\emptyset.\nonumber
\end{equation}

Now we will show that this intersection is transversal. First note that from
the analyticity of the functions $v^{s}(t_{0},e)$ and $v^{u}(t_{0},e)$ using
the same argument as in the proof of Lemma \ref{lem:Melnikov} we also have
that
\begin{equation}
\frac{\partial}{\partial t}\left(  h^{u}(t_{0},e)-h^{s}(t_{0},e)\right)
=e\frac{\partial}{\partial t}M_{\mu_{k}}(t)+O(e\mu_{k}^{-1/3}R(\mu_{k}))+o(e).
\label{eq:partial-d}%
\end{equation}
We know that prior to our perturbation $W_{l(r)}^{s}$ and $W_{l(r)}^{u}$
intersect transversally at $v^{0}(t_{0})$ (\ref{eq:voto}). This intersection
is not transversal in the full extended phase space, it is only transversal in
the constant energy manifold
\begin{equation}
M=\{(\bar{x},\bar{y},\bar{p}_{x},\bar{H},t)|\bar{H}=\bar{H}(\mu_{k}%
,l(r))\}\subset\mathbb{R}^{3}\times\{\bar{H}(\mu_{k},l(r))\}\times
S^{1}.\nonumber
\end{equation}
To be more precise, we know that $v^{0}(t_{0})\in W_{l_{t_{0}}^{0}(r)}^{s}\cap
W_{l_{t_{0}}^{0}(r)}^{u}$, that $W_{l_{t_{0}}^{0}(r)}^{s},$ $W_{l_{t_{0}}%
^{0}(r)}^{u}\subset\Sigma_{t_{0}}$ and that \cite{Simo}
\begin{equation}
T_{v^{0}(t_{0})}(W_{l_{t_{0}}^{0}(r)}^{s})+T_{v^{0}(t_{0})}(W_{l_{t_{0}}%
^{0}(r)}^{u})=\mathbb{R}^{3}\times\{0\}\times\{0\}.\nonumber
\end{equation}
These properties are preserved under small perturbation $e>0$, hence for
$v^{s}(t_{0},e)=v^{u}(t_{0},e)=:v(t_{0},e)$ we have%
\begin{equation}
T_{v(t_{0},e)}(W_{l_{t_{0}}^{e}(r)}^{s})+T_{v(t_{0},e)}(W_{l_{t_{0}}^{e}%
(r)}^{u})=\mathbb{R}^{3}\times\{0\}\times\{0\}.\nonumber
\end{equation}
We need to show that we also have transversality on the $t_{0}$ and energy
coordinate. For a fixed $e$ the curves $v^{s}(t,e)$ and $v^{u}(t,e)$ belong to
$W_{l^{e}(r)}^{s}$ and $W_{l^{e}(r)}^{u}$ respectively. At the time $t=t_{0}$
for which $v^{s}(t_{0},e)=v^{u}(t_{0},e)=v(t_{0},e)$ we have $M_{\mu_{k}%
}^{\prime}(t_{0})\neq0.$ This means that using (\ref{eq:partial-d}), for
sufficiently small $e,$%
\begin{align*}
\frac{\partial}{\partial t}\left(  \pi_{H}(v^{s}(t,e))-\pi_{H}(v^{u}%
(t,e))\right)  |_{t=t_{0}}  &  =\frac{\partial}{\partial t}\left(
h^{u}(t,e)-h^{s}(t,e)\right)  |_{t=t_{0}}\\
&  =eM_{\mu_{k}}^{\prime}(t_{0})+O(e\mu_{k}^{-1/3}R(\mu_{k}))+o(e)\\
&  \neq0.
\end{align*}
We also have for $i\in\{u,s\},$ $\frac{\partial}{\partial t}\left(  \pi
_{t_{0}}v^{i}(t,e)\right)  =\frac{\partial}{\partial t}t=1.$ This since
$\frac{d}{dt}v^{s}(t,e)|_{t=t_{0}}\in T_{v(t_{0},e)}(W_{l^{e}(r)}^{s})$ and
$\frac{d}{dt}v^{u}(t,e)|_{t=t_{0}}\in T_{v(t_{0},e)}(W_{l^{e}(r)}^{u})$
implies transversality, which finishes our proof.

The order of choice of parameters in the above argument is important, so let
us quickly run through how it should be conducted. We first choose
sufficiently small $\mu_{k}$ so that we can apply Theorem \ref{th:KAM for L2}.
Then by choosing small $R(\mu_{k})$ we ensure that $\mu_{k}^{-1/3}R(\mu_{k})$
is sufficiently small compared with $M_{\mu_{k}}$ and $M_{\mu_{k}}^{\prime}$.
We then choose $e$ so that $e\mu_{k}^{-1/3}$ is sufficiently small so that we
have transversal intersections of $\pi_{x,p_{x}}(W_{l^{e}(r)}^{u}\cap
\Sigma_{t_{0}}\cap\{\bar{y}=0\})$ and $\pi_{x,p_{x}}(W_{l^{e}(r)}^{s}%
\cap\Sigma_{t_{0}}\cap\{\bar{y}=0\})$. The parameter $e$ needs also to be
sufficiently small so that $M_{\mu_{k}}(t_{0})$ and $M_{\mu_{k}}^{\prime
}(t_{0})$ dominate in (\ref{eq:H-difference}) and (\ref{eq:partial-d})
respectively. \begin{figure}[ptb]
\begin{center}
\includegraphics[
height=1.3in
]{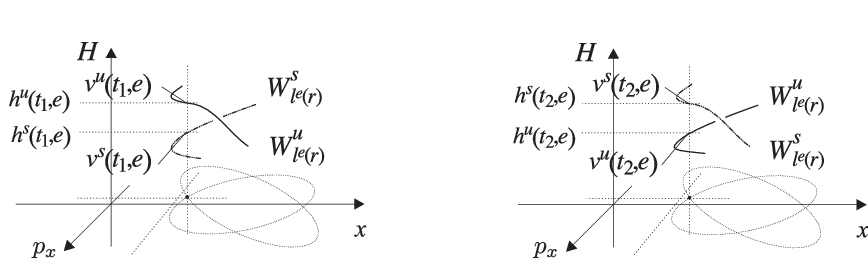}
\end{center}
\caption{The manifolds $W_{l^{e}(r)}^{s}$ and $W_{l^{e}(r)}^{u}$ intersected
with $\Sigma_{t_{1}},\Sigma_{t_{2}}$ and $\{y=0\},$ in the $x,p_{x},H$
coordinates for $t_{1}<t_{0}<t_{2}.$}%
\label{fig:p(c,e)-3d}%
\end{figure}
\end{proof}

\begin{remark}
\label{rem:size-of-R}In the proof of Theorem \ref{th:Melnikov} we
see that we need to choose the radius $R(\mu_{k})$ to be
sufficiently small so that $\mu_{k}^{-1/3}R(\mu_{k})$ is small in
comparison to the Melnikov integral $M_{\mu_{k}}$ and its
derivative. Since we do not have uniform bounds on the size of the
Melnikov integral with respect to $\mu_{k},$ from our argument we
cannot say how small $R(\mu_{k})$ needs to be. From our numerical
investigation which will follow in Table 2 we can see that for
$t_{0}$ for which we have a simple zero of the Melnikov integral,
the bound on the derivative is independent of $\mu_{k}$. This
means that we need to choose the radius $R(\mu_{k})\mu_{k}^{1/3}$
to be sufficiently small (hence $R(\mu_{k})$ converges to zero
with $\mu_{k}$ going to zero).
\end{remark}

\begin{corollary}
\label{cor:melnikov}If $\mu_{k}$ is sufficiently small and the Melnikov
integral has a simple zero then for sufficiently small $\mu_{k}^{-1/3}%
R(\mu_{k})$ there exists a $\zeta>0$, such that for any two radii $r_{1}$ and
$r_{2}$ from $\mathfrak{C}\cap\lbrack R,R(\mu_{k})]$%
\begin{equation}
|r_{1}-r_{2}|<\zeta e\mu_{k}^{-1/3},\nonumber
\end{equation}
the manifolds $W_{l^{e}(r_{i})}^{s}$ and $W_{l^{e}(r_{j})}^{u}$ intersect
transversally for $i,j\in\{1,2\}$.
\end{corollary}

\begin{proof}
The proof of this fact is a mirror argument to the proof of Theorem
\ref{th:Melnikov}. Below we restrict our attention to pointing out the
difference we have connected with the derivation of (\ref{eq:H-difference}) in
the setting where we have two radii.

Let $h_{1}^{u}(t_{0},e)$ and $h_{2}^{s}(t_{0},e)$ stand for energies of points
of potential intersection (constructed analogously to $h^{u}(t_{0},e)$ and
$h^{s}(t_{0},e)$ in (\ref{eq:hu-hs})). Let $r_{1}^{u}$ and $r_{2}^{s}$ stand
for radii constructed analogously to $r^{u}$ and $r^{s}$ (see (\ref{eq:ru-rs}%
)), but coming from the unstable and stable manifold of $l^{e}(r_{1})$ and
$l^{e}(r_{2})$ respectively. Using Lemma \ref{lem:dist(Lc,L2)} and the fact
that $r_{1}^{u}-r_{1}=O(e\mu_{k}^{-1/3})$ and $r_{2}^{s}-r_{2}=O(e\mu
_{k}^{-1/3})$ we have
\begin{align*}
|\bar{H}(\mu_{k},l(r_{1}^{u}))-\bar{H}(\mu_{k},l(r_{2}^{s}))|  &  \leq
O(|(r_{1}^{u})^{2}-(r_{2}^{s})^{2}|)\\
&  \leq O(R(\mu_{k})|r_{1}^{u}-r_{2}^{s}|)\\
&  \leq O(R(\mu_{k})|r_{1}-r_{2}|)+O(e\mu_{k}^{-1/3}R(\mu_{k})).
\end{align*}
Using an identical argument to the derivation of (\ref{eq:H-difference}) this
gives us
\begin{equation}
h_{1}^{u}(t_{0},e)-h_{2}^{s}(t_{0},e)=eM_{\mu_{k}}(t_{0})+O(R(\mu_{k}%
)|r_{1}-r_{2}|)+O(e\mu_{k}^{-1/3}R(\mu_{k}))+o(e).\nonumber
\end{equation}
Using this estimate and following the proof of Theorem \ref{th:Melnikov} we
obtain our claim.
\end{proof}

\begin{remark}
\label{rem:melnikov} Mirror arguments to the proof of Theorem
\ref{th:Melnikov} and Corollary \ref{cor:melnikov} give transversal
intersections of invariant manifolds for invariant tori of the PRE3BP with
$\mu=\mu_{k}$ for $k=2,3,\ldots$. Here we state this as a separate remark
since Theorem \ref{th:Melnikov} and Corollary \ref{cor:melnikov} are fully
rigorous and do not rely on any numerical computations. For the argument with
an arbitrary $\mu_{k}$ we would need to use the fact that in the PRC3BP we
have a twist property on the family of Lyapounov orbits and apply Remark
\ref{rem:KAM-for-all-k}. The twist for arbitrary $\mu_{k}$ has only been
demonstrated numerically (see Table 1).
\end{remark}

\section{Computation of the Melnikov integral.}

\label{sec:meln-comp}

In this section we will demonstrate that for $t_{0}=0$ and for all parameters
$\mu_{k}$ from Theorem \ref{lem homoclinic orbit for uk} the Melnikov integral
$M_{\mu_{k}}(t_{0})$ (\ref{eq:Melnikov}) is zero and also that $\frac
{dM_{\mu_{k}}}{dt_{0}}(0)\neq0.$ The fact that the Melnikov integral is zero
will follow directly from the $S$-symmetry (\ref{eq:S-sym}) of the homoclinic
orbit $q_{\mu_{k}}^{0}(t)$. The fact that $\frac{dM_{\mu_{k}}}{dt_{0}}%
(0)\neq0$ will be demonstrated numerically. We will compute the integral for
the first few parameters $\mu_{k}$ and then demonstrate that for sufficiently
small parameters $\mu_{k}$ the integral converges to an integral along an
unstable manifold of the Hill's problem.

\subsection{The Melnikov integral and its derivative at $t_{0}=0$}

We start with a lemma which ensures the convergence of the Melnikov integral
(\ref{eq:Melnikov}).

\begin{lemma}
\label{lem:Meln-potential}The Melnikov integral (\ref{eq:Melnikov}) and its
derivative is absolutely convergent uniformly with respect to $t_{0}$. The
Melnikov function can be expressed as
\begin{equation}
M_{\mu_{k}}\left(  t_{0}\right)  =\int_{-\infty}^{+\infty}\left[
\frac{\partial\bar{G}}{\partial t}\left(  \mu_{k},\bar{q}_{\mu_{k}}^{0}\left(
t\right)  ,t+t_{0}\right)  -\frac{\partial\bar{G}}{\partial t}\left(  \mu
_{k},\bar{L}_{2}^{\mu_{k}},t+t_{0}\right)  \right]  dt,
\label{eq:M-from-potential}%
\end{equation}
and also%
\begin{equation}
\frac{dM_{\mu_{k}}}{dt}\left(  t_{0}\right)  =\int_{-\infty}^{+\infty}\left[
\frac{\partial^{2}\bar{G}}{\partial t^{2}}\left(  \mu_{k},\bar{q}_{\mu_{k}%
}^{0}\left(  t\right)  ,t+t_{0}\right)  -\frac{\partial^{2}\bar{G}}{\partial
t^{2}}\left(  \mu_{k},\bar{L}_{2}^{\mu_{k}},t+t_{0}\right)  \right]  dt.
\label{eq:der-M}%
\end{equation}

\end{lemma}

\begin{proof}
The orbit $\bar{q}_{\mu_{k}}^{0}\left(  t\right)  $ is the homoclinic orbit to
the Libration point $\bar{L}_{2}^{\mu_{k}}.$ Let us note that the velocity
$\bar{x}^{\prime}$ and $\bar{y}^{\prime}$ of $\bar{q}_{\mu_{k}}^{0}(t)$
exponentially tends to zero as $t$ tends to plus infinity and minus infinity.
Moreover the partial derivatives of $\bar{G}$ on $\bar{q}_{\mu_{k}}^{0}(t)$
are uniformly bounded. This means that the integral over
\begin{equation}
\int_{-\infty}^{+\infty}|\{\bar{H},\bar{G}\}|(\bar{q}_{\mu_{k}}^{0}%
(t),t+t_{0})dt=\int_{-\infty}^{+\infty}|\bar{x}^{\prime}\frac{\partial\bar{G}%
}{\partial\bar{x}}+\bar{y}^{\prime}\frac{\partial\bar{G}}{\partial\bar{y}%
}|(\mu_{k},\bar{q}_{\mu_{k}}^{0}(t),t+t_{0})dt,\nonumber
\end{equation}
is convergent uniformly with respect to $t_{0}$.

The orbit $\bar{q}_{\mu_{k}}^{0}\left(  t\right)  $ is the solution of the
PRC3BP, hence differentiating gives
\begin{equation}
\frac{d\bar{G}}{dt}\left(  \mu_{k},\bar{q}_{\mu_{k}}^{0}\left(  t\right)
,t+t_{0}\right)  =\frac{\partial\bar{G}}{\partial t}\left(  \mu_{k},\bar
{q}_{\mu_{k}}^{0}\left(  t\right)  ,t+t_{0}\right)  +\left\{  \bar{G},\bar
{H}\right\}  \left(  \mu_{k},\bar{q}_{\mu_{k}}^{0}\left(  t\right)
,t+t_{0}\right)  .\label{eq:derivative-poisson}%
\end{equation}
From (\ref{eq:derivative-poisson}) we have
\begin{align}
M_{\mu_{k}}(t_{0})= &  \int_{-\infty}^{+\infty}\left\{  \bar{H},\bar
{G}\right\}  (\mu_{k},\bar{q}_{\mu_{k}}^{0}(t),t+t_{0})dt\nonumber\\
= &  \lim_{T\rightarrow\infty}\int_{-T}^{T}\left(  \frac{\partial\bar{G}%
}{\partial t}\left(  \mu_{k},\bar{q}_{\mu_{k}}^{0}\left(  t\right)
,t+t_{0}\right)  -\frac{d\bar{G}}{dt}\left(  \mu_{k},\bar{q}_{\mu_{k}}%
^{0}\left(  t\right)  ,t+t_{0}\right)  \right)  dt\label{eq:M(0)-form1}\\
= &  \lim_{T\rightarrow\infty}\left[  \bar{G}\left(  \mu_{k},\bar{q}_{\mu_{k}%
}^{0}\left(  -T\right)  ,-T+t_{0}\right)  -\bar{G}\left(  \mu_{k},\bar{q}%
_{\mu_{k}}^{0}\left(  T\right)  ,T+t_{0}\right)  \right.  \nonumber\\
&  \left.  +\int_{-T}^{T}\frac{\partial\bar{G}}{\partial t}\left(  \mu
_{k},\bar{q}_{\mu_{k}}^{0}\left(  t\right)  ,t+t_{0}\right)  dt\right]
.\nonumber
\end{align}
To complete the proof of (\ref{eq:M-from-potential}) observe that that from
$\lim_{T\rightarrow\pm\infty}\bar{q}_{\mu_{k}}^{0}\left(  T\right)  =\bar
{L}_{2}^{\mu_{k}}$ it follows that
\begin{eqnarray}
\lim_{T\rightarrow\infty}\left(  \int_{-T}^{T}\frac{\partial\bar{G}}{\partial
t}\left(  \mu_{k},\bar{L}_{2}^{\mu_{k}},t+t_{0}\right)  dt-\bar{G}\left(
\mu_{k},\bar{q}_{\mu_{k}}^{0}\left(  T\right)  ,T+t_{0}\right)  \right.
\label{eq:M(0)-form2}\\
\left.  \phantom{\int_T^T}+\bar{G}\left(  \mu_{k},\bar{q}_{\mu_{k}}^{0}\left(
-T\right)  ,-T+t_{0}\right)  \right)  =0\nonumber
\end{eqnarray}
uniformly with respect to $t_{0}$.

From (\ref{eq:M(0)-form1}) and (\ref{eq:M(0)-form2}) we obtain
(\ref{eq:M-from-potential}).

To prove (\ref{eq:der-M}) it is enough to observe that the formal
integration of (\ref{eq:M-from-potential}) is correct, because the
integral on the right hand side of (\ref{eq:der-M}) is uniformly
convergent with respect to $t_{0}$ for the same reasons as the
integral in formula (\ref{eq:M-from-potential}).
\end{proof}

It turns out that the computation of the Melnikov integral is not
the major obstacle. The fact that we have a zero for $t_{0}=0$
follows directly from the $S$-symmetry (\ref{eq:S-sym}) of the
homoclinic orbit $q_{\mu_{k}}^{0}(t)$ and $G$. This fact is shown
in the below lemma . Later though we
will need to show that this zero is nontrivial by computing
$\frac{dM_{\mu_{k}}}{dt}(0)$, which turns out to be a much harder
task. The lemma also provides the formula for the needed integral.

\begin{lemma}
\label{lem:Meln1-symmetry}The Melnikov integral (\ref{eq:Melnikov}) at
$t_{0}=0$ is equal to zero and
\begin{equation}
\frac{dM_{\mu_{k}}}{dt}\left(  0\right)  =-2\mu_{k}^{-2/3}\int_{-\infty}%
^{0}\left(  G\left(  \mu_{k},q_{\mu_{k}}^{0}\left(  t\right)  ,t\right)
-G\left(  \mu_{k},L_{2}^{\mu_{k}},t\right)  \right)  dt.
\label{eq:dM/dt=2*int}%
\end{equation}

\end{lemma}

\begin{proof}
First let us observe that from (\ref{eq:G-bar})
\begin{equation}
\bar{G}\left(  \mu_{k},\bar{q}_{\mu_{k}}^{0}\left(  t\right)  ,t\right)
=\mu_{k}^{-2/3}G\left(  \mu_{k},q_{\mu_{k}}^{0}\left(  t\right)  ,t\right)
.\nonumber
\end{equation}
The orbit $q_{\mu_{k}}^{0}\left(  t\right)  $ and fixed point $L_{2}^{\mu_{k}%
}$ are $S$-symmetric with respect to the symmetry (\ref{eq:S-sym}). From
(\ref{eq:S-def}), (\ref{eq:Gdef}), by direct computation one can check that
$\frac{\partial G}{\partial t}\left(  \mu,S\left(  \cdot\right)  ,-t\right)
=-\frac{\partial G}{\partial t}\left(  \mu,\cdot,t\right)  $ and that
$\frac{\partial^{2}G}{\partial t^{2}}\left(  \mu,S\left(  \cdot\right)
,-t\right)  =\frac{\partial^{2}G}{\partial t^{2}}\left(  \mu,\cdot,t\right)
\ $hence we have
\begin{align*}
\int_{-\infty}^{0}  &  \left(  \frac{\partial\bar{G}}{\partial t}\left(
\mu_{k},\bar{q}_{\mu_{k}}^{0}\left(  t\right)  ,t\right)  -\frac{\partial
\bar{G}}{\partial t}\left(  \mu_{k},\bar{L}_{2}^{\mu_{k}},t\right)  \right)
dt=\\
&  =\mu_{k}^{-2/3}\int_{0}^{+\infty}\left(  \frac{\partial G}{\partial t}%
(\mu_{k},q_{\mu_{k}}^{0}(-t),-t)-\frac{\partial G}{\partial t}(\mu_{k}%
,L_{2}^{\mu_{k}},-t)\right)  dt\\
&  =\mu_{k}^{-2/3}\int_{0}^{+\infty}\left(  \frac{\partial G}{\partial
t}\left(  \mu_{k},S\left(  q_{\mu_{k}}^{0}\left(  t\right)  \right)
,-t\right)  -\frac{\partial G}{\partial t}(\mu_{k},S\left(  L_{2}^{\mu_{k}%
}\right)  ,-t)\right)  dt\\
&  =-\mu_{k}^{-2/3}\int_{0}^{+\infty}\left(  \frac{\partial G}{\partial
t}\left(  \mu_{k},q_{\mu_{k}}^{0}\left(  t\right)  ,t\right)  -\frac{\partial
G}{\partial t}\left(  \mu_{k},L_{2}^{\mu_{k}},t\right)  \right)  dt,
\end{align*}
which gives $M_{\mu_{k}}(0)=0.$ From an analogous computation using
$\frac{\partial^{2}G}{\partial t^{2}}\left(  S\left(  \cdot\right)
,-t\right)  =\frac{\partial^{2}G}{\partial t^{2}}\left(  \cdot,t\right)  $
follows%
\begin{equation}
\frac{dM_{\mu_{k}}}{dt}\left(  0\right)  =2\mu_{k}^{-2/3}\int_{-\infty}%
^{0}\left(  \frac{\partial^{2}G}{\partial t^{2}}\left(  \mu_{k},q_{\mu_{k}%
}^{0}\left(  t\right)  ,t\right)  -\frac{\partial^{2}G}{\partial t^{2}}\left(
\mu_{k},L_{2}^{\mu_{k}},t\right)  \right)  dt.\nonumber
\end{equation}
Form (\ref{eq:Gdef}) we have $\frac{\partial^{2}G}{\partial t^{2}}=-G,$ which
gives (\ref{eq:dM/dt=2*int}).
\end{proof}

\begin{remark}\label{rem:rig-num}
The verification of the fact that $\frac{dM_{\mu_{k}}}{dt}(0)$ is
nonzero is not straightforward. In this paper we will restrict
ourselves to numerical verification of this fact. We would like to
highlight that for a given parameter $\mu_{k}$ it is possible to
obtain a rigorous-computer-assisted estimate on
$\frac{dM_{\mu_{k}}}{dt}(0)$. To do so one first needs to obtain a
rigorous bound $[$\underline{$\mu$}$_{k},\overline{\mu}_{k}]$
which contains the parameter $\mu_{k}$ for which we have a
homoclinic orbit. Then one needs to obtain rigorous enclosures on
the trajectories $q_{\mu}^{0}(t)$ for all
$\mu\in\lbrack$\underline{$\mu$}$_{k},\overline{\mu}_{k}]$. Using
these, a bound on $\frac{dM_{\mu_{k}}}{dt}(0)$ can be computed
from (\ref{eq:dM/dt=2*int}).

We have successfully conducted such computations for the parameter
$\mu_{2}$. We have used the fact that if one extends the system by
including $\mu$ as an additional variable, then for any interval
$I$ the set $\{L_{2}^{\mu}\}_{\mu\in I}$ is a normally hyperbolic
invariant manifold. We have applied a topological method, given in
\cite{cmfld}, \cite{nhim}, for detection of normally hyperbolic
manifolds, combined with a parameterization method
\cite{llave-parameter}. Based on these we proved the following. We
have shown that the parameter $\mu_{2}$ for which we have the
homoclinic orbit is contained in $\left[
0.0042538631,0.0042538639\right]$. We then proved that
$\frac{dM_{\mu_{2}}}{dt}(0)\in\left[
1.301020122,1.865308899\right]$. A detailed proof of this fact,
along with results for other parameters, will be the subject of a
forthcoming publication.
\end{remark}

In Table 2 we enclose the (nonrigorous) numerical results for the computation
of (\ref{eq:dM/dt=2*int}) obtained for $k$ up to $13.\mathbf{\ }$

\begin{samepage}
\begin{center}
\begin{tabular}[c]{l | r  }
$k$ & $\frac{dM_{\mu_k}}{dt}(0)$\\
&  \\ \hline
2   & 1.57396 \\
3   &  -0.396727\\
4   & 0.395119 \\
5   &  -0.396931\\
6   & 0.395784 \\
7   & -0.395511 \\
8   & 0.393389 \\
9   & -0.393253 \\
10 & 0.390924 \\
11 & -0.391194 \\
12 &   0.388961\\
13 & -0.389459
\end{tabular}

\textbf{Table 2.} Numerical results for the derivation of
$\frac{dM_{\mu_k}}{dt}(0)$ for various mass parameters. \end{center}
\end{samepage}

To obtain the numerical results from Table 2 we have used a
parameterization method \cite{llave-parameter} to obtain an
expansion of the manifold around $L_{2}^{\mu_{k}}$ as a polynomial
of degree 20. Then we integrated the system numerically using a
Taylor method of order 20. Numerical evidence points to
$M_{\mu_{k}}$ having a nontrivial zero for $t=0$.

\section{Transition Chains and Main Result}

\label{sec:tranchains}

In this section we will show that there exists a sequence of Lyapounov orbits
$l(r_{i})$ for $i=1,\ldots,N$ which survive for a sufficiently small
perturbation $e,$ and such that their stable and unstable manifolds intersect
transversally%
\begin{equation}
W_{l^{e}(r_{i})}^{s}\pitchfork W_{l^{e}(r_{i+1})}^{u}\quad\text{and}\quad
W_{l^{e}(r_{i})}^{u}\pitchfork W_{l^{e}(r_{i+1})}^{s}.
\label{eq:trans-chain-inter}%
\end{equation}
Such sequences are referred to as transition chains. Existence of such chains
ensures that we have an orbit which shadows the homo- or heteroclinic
connections between the surviving tori generating rich symbolic dynamics (see
\cite{Gidea-Llave},\cite{Gidea-Robinson},\cite{Gidea}). Apart from this, from
our argument it will also follow that for any surviving orbit $l(r)$ we also
have transversal intersection of its stable and unstable manifold%
\begin{equation}
W_{l^{e}(r)}^{u}\pitchfork W_{l^{e}(r)}^{s}. \label{eq:trans-chain-inter-2}%
\end{equation}
This ensures that the chaotic dynamics of the PRC3BP which is implied by the
existence of a transversal homoclinic orbit to $l(r)$ (see Theorem
\ref{lem: transversality-Simo} and Remark \ref{rem:symb-dyn}) survives.

We are now ready to rigorously reformulate our main theorem (Theorem
\ref{thm:main}).

\begin{theorem}
\label{th:Main-Detail}For any $\mu$ from the sequence of masses $\{\mu
_{k}\}_{k=2}^{\infty}$ from Theorem \ref{lem homoclinic orbit for uk}, if the
twist property is satisfied, then there exists a radius $R_{\text{Hill}}$,
which is independent of $\mu$, such that for eccentricities $e$ of the
elliptic problem, with $e\mu^{-2/3}<\kappa$ (see \ref{eq:kappa-e-bound} for
interpretation of $\kappa$) and sufficiently small $e\mu^{-1/3}$, there exists
a Cantor set $\mathfrak{C}\subset\lbrack0,R_{\text{Hill}}]$ such that for all
$r\in\mathfrak{C}$ the Lyapounov orbits $l(r)$ are perturbed to invariant tori
$l^{e}(r)$ in the extended phase space. Moreover, if the derivative of the
Melnikov integral (\ref{eq:dM/dt=2*int}) is nonzero, then there exists a
radius $R(\mu)$ of order at most $O(\mu^{1/3})$ and a transition chain
$l^{e}(r_{i})$ for a sequence of radii $0<r_{1}<r_{2}<\ldots<r_{N}<R(\mu)$
such that the difference  $r_{N}-r_{1}$ is of order $\left(  e\mu
^{-1/3}\right)  ^{1/2}$.

The transversal intersections of the transition chain
(\ref{eq:trans-chain-inter}), (\ref{eq:trans-chain-inter-2}) lead to a homo-
and heteroclinic tangle of the stable/unstable manifolds of $l^{e}(r_{i} )$,
which in turn leads to symbolic dynamics involving diffusion in energy.
\end{theorem}

\begin{remark}
Let us note now that we can verify assumptions of the above theorem based on
some numerical results. The twist property for the family of Lyapounov orbits
has been rigorously proved for sufficiently small $\mu_{k}$ in Theorem
\ref{th:KAM for L2}, but the fact that we have twist for $\mu_{k}$ with
$k=1,2,\ldots$ has only been demonstrated numerically (see Section
\ref{sec:twist-L2}, Table 1 and Remark \ref{rem:KAM-for-all-k}). Secondly,
assumptions of Theorem \ref{th:Main-Detail} require that the derivative of the
Melnikov function (\ref{eq:dM/dt=2*int}) at zero is nonzero. This has only
been demonstrated numerically in Section \ref{sec:meln-comp}.

We believe that the above can be verified using rigorous-computer-assisted
methods. This is currently a subject of ongoing work (see Remark
\ref{rem:rig-num}).
\end{remark}

\begin{remark}
All radiuses considered in Theorem \ref{th:Main-Detail} are given in the
Hill's coordinates (\ref{eq:Hill-coordinates}). This means that in the
original coordinates of our system (given by the Hamiltonian
(\ref{eq:H-PRE3BP})) the radiuses are reduced by a factor of $\mu^{1/3}$.
\end{remark}

\begin{proof}
[Proof of Theorem \ref{th:Main-Detail}]Let us fix a $\mu=\mu_{k}.$ For
sufficiently small $\mu$ by Theorem \ref{th:KAM for L2}, and if the twist
condition holds for all $\mu_{k}$ by Remark \ref{rem:KAM-for-all-k}, we have
the radius $R_{\text{Hill}}$ and a Cantor $\mathfrak{C}\subset\lbrack
0,R_{\text{Hill}}]$ of radii for which the Lyapounov orbits survive the
perturbation. From Theorem \ref{th:Melnikov}, Remark \ref{rem:melnikov} (see
also Corollary \ref{cor:melnikov}) it follows that by choosing $R(\mu
)<R_{\text{Hill}}$, for which $R(\mu)\mu^{-1/3}$ is sufficiently small, we
know that there exists a $\zeta>0$ such that for radii $r_{1},r_{2}%
\in\mathfrak{C}$ such that $|r_{1}-r_{2}|<\zeta e\mu^{-1/3}$, for $e$ with
sufficiently small $e\mu^{-1/3},$ we have%
\begin{equation}
W_{l^{e}(r_{i})}^{u}\pitchfork W_{l^{e}(r_{i})}^{s}\quad\text{for}\quad
i\in\{1,2\}.\nonumber
\end{equation}
We now need to show that we can find a sequence $r_{1}<r_{2}<\ldots<r_{N}$
such that $r_{i}\in\mathfrak{C}$ and the difference $r_{N}-r_{1}$ is of order
$\left(  e\mu^{-1/3}\right)  ^{1/2},$ for which the gaps between $r_{i}$ and
$r_{i+1}$ are smaller than $\zeta e\mu^{-1/3}.$ The existence of such a
sequence follows from Proposition \ref{prop:KAM-gaps}.

The last claim of Theorem \ref{th:Main-Detail} follows from \cite{Gidea-Llave}%
, \cite{Gidea-Robinson}.
\end{proof}

\section{Concluding remarks, future work}

In this paper we have shown that the chaotic dynamics observed for the planar
restricted circular three body problem survives the perturbation into the
planar restricted elliptic three body problem, when its eccentricity is
sufficiently small. We have also shown that this dynamics is extended to
include diffusion in energy. The diffusion proved in this paper covers a small
range of energies. This is due to the fact that in our argument we use a
Melnikov type method which does not allow us to jump between the "large gaps"
between the KAM tori. An interesting problem which could be addressed is
whether these large gaps can be overcome (this potentially could be done using
techniques similar to \cite{large-gap} or \cite{Gidea}).

Our result holds only for a specific family $\{\mu_{k}\}$ of masses of the
primaries. The choice of these masses is such that they ensure the existence
of the homoclinic orbit to $L_{2}^{\mu_{k}}$, which is then used for the
Melnikov argument. An interesting question is whether one can observe similar
dynamics in real life setting, say in the Jupiter-Sun system. In such a case
we will no longer have a homoclinic connection for the point $L_{2}^{\mu}$.
For the (circular) Jupiter-Sun system though we know that we have a
transversal homoclinic connection for Lyapounov orbits (see \cite{Marsden} and
\cite{Zgliczyn-Wilczak}). Such orbits could be used for a similar
construction. Our argument also required that we have sufficiently small
eccentricities. It would be interesting to find out if the dynamics persists
for the actual eccentricity of the Jupiter-Sun system. For this problem it is
quite likely that applying the mechanism discussed in this paper would be very
hard. Our argument relies on the use of the KAM theorem, which works for
sufficiently small perturbations. To apply it for an explicit eccentricity
seems a difficult task. Other methods could be exploited though. Instead of
proving the persistence of the tori and trying to detect intersections of
their invariant manifolds, one could focus on detection of symbolic dynamics
for the diffusing orbits in the spirit of \cite{Zgliczyn-Wilczak}. This seems
a far more realistic target for the near future and is being currently
considered as an extension of this work.

\section{Acknowledgements}

The authors would like to thank the anonymous reviewers for their valuable
comments and suggestions which helped them improve the quality of the paper.

\section{Appendix}

\subsection{Splitting of manifolds associated to Lyapounov orbits of the
PRC3BP}

Here we investigate the dependence on parameter $\mu$ of the splitting of the
intersection of curves from Theorem \ref{lem: transversality-Simo}. Let
$\mu=\mu_{k},$ which means that in (\ref{eq:Mf-implicit}) $\alpha=0.$ The
transversal intersection of the curves is on $\{\dot{x}=0\}$. The curve
intersects $\{\dot{x}=0\}$ for $\sigma=\sigma_{0}$ for which
\begin{equation}
(K_{1}\cos\tau\cos\sigma_{0}-K_{2}\sin\tau\sin\sigma_{0})=O(\mu^{4/3}%
).\label{eq:sig0-implicit}%
\end{equation}
This is because from (\ref{eq:Mf-implicit}) only then can we have
$M_{f}(\sigma_{0})$ sufficiently close to zero so that $\dot{x}=0.$ In
particular, (\ref{eq:Mf-implicit}) and (\ref{eq:sig0-implicit}) gives
\begin{equation}
M_{f}(\sigma_{0})=O(\mu^{2/3}).\label{eq:Mf-sig0}%
\end{equation}
By implicit function theorem we know that
\begin{align}
\gamma(\mu,\sqrt{\Delta C}) &  :=\frac{\partial M}{\partial\sigma}\left(
\sigma_{0}\right)  \label{eq:gamma-est-mu}\\
&  =-\frac{\mu^{-2/3}\sqrt{\Delta C}3M\left(  N+2M\cos\tau\right)
^{-1}\left(  -K_{1}\cos\tau\sin\sigma_{0}-K_{2}\sin\tau\cos\sigma_{0}\right)
}{N+2M\cos(M_{f}(\sigma_{0}))}\nonumber\\
&  \in\mu^{-2/3}\sqrt{\Delta C}\left[  a_{1},a_{2}\right]  ,\nonumber
\end{align}
where $a_{1},a_{2}$ are constants, $0\notin\left[  a_{1},a_{2}\right]  $ (for
sufficiently small $\mu$ and $\Delta C$ these constants are arbitrarily close
to one another. For their precise values one would need to substitute the
constants $M,$ $\tau$, $K_{1},$ $K_{2},$ $N$ from \cite{Simo} into
(\ref{eq:gamma-est-mu})). By (\ref{eq:sig0-implicit}), (\ref{eq:Mf-sig0}),
(\ref{eq:gamma-est-mu})
\begin{align*}
&  \frac{\partial x}{\partial\sigma}\left(  \sigma_{0}\right)  \\
&  =\sqrt{\Delta C}\left(  N+2M\cos\tau\right)  ^{-1}N\left(  -\sin\left(
M_{f}(\sigma_{0})\right)  \gamma(\mu,\sqrt{\Delta C})\right)  \\
&  \cdot(K_{1}\cos\tau\cos\sigma_{0}-K_{2}\sin\tau\sin\sigma_{0})\\
&  +\sqrt{\Delta C}\left(  N+2M\cos\tau\right)  ^{-1}\left(  2M+N\cos
M_{f}(\sigma_{0})\right)  \left(  -K_{1}\cos\tau\sin\sigma_{0}-K_{2}\sin
\tau\cos\sigma_{0}\right)  \\
&  +\mu^{1/3}M\sin\left(  M_{f}(\sigma_{0})\right)  \gamma(\mu,\sqrt{\Delta
C})\\
&  +\mu^{2/3}\left\{  -\frac{2MN}{3}\sin\left(  M_{f}\left(  \sigma
_{0}\right)  \right)  \gamma(\mu,\sqrt{\Delta C})+M^{2}2\sin\left(
M_{f}\left(  \sigma_{0}\right)  \right)  \cos\left(  M_{f}\left(  \sigma
_{0}\right)  \right)  \gamma(\mu,\sqrt{\Delta C})\right\}  \\
&  \in\sqrt{\Delta C}\left[  b_{1},b_{2}\right]  ,
\end{align*}
where $b_{1},b_{2}$ are constants, $0\notin\left[  b_{1},b_{2}\right]  $ (the
second term in the above equation plays a dominant role). Also%
\begin{align*}
\frac{\partial\dot{x}}{\partial\sigma}\left(  \sigma_{0}\right)   &
=\cos\left(  M_{f}(\sigma_{0})\right)  \gamma(\mu,\sqrt{\Delta C})\left[
\frac{\sqrt{\Delta C}N\left(  K_{1}\cos\tau\cos\sigma_{0}-K_{2}\sin\tau
\sin\sigma_{0}\right)  }{N+2M\cos\tau}\right.  \\
&  \left.  +\mu^{1/3}M+\mu^{2/3}\left\{  \frac{MN}{3}+2M^{2}\cos M_{f}%
+\frac{M}{3}\alpha\right\}  \right]  \\
&  +\sin\left(  M_{f}(\sigma_{0})\right)  [\sqrt{\Delta C}N\left(
N+2M\cos\tau\right)  ^{-1}\left(  -K_{1}\cos\tau\sin\sigma_{0}-K_{2}\sin
\tau\cos\sigma_{0}\right)  \\
&  -2M^{2}\sin\left(  M_{f}(\sigma_{0})\right)  \gamma(\mu,\sqrt{\Delta C})]\\
&  \in\mu^{-1/3}\sqrt{\Delta C}[c_{1},c_{2}],
\end{align*}
where $c_{1},c_{2}$ are constants, $0\notin\left[  c_{1},c_{2}\right]  .$ This
means that
\begin{equation}
\frac{\partial x}{\partial\sigma}(\sigma_{0})/\frac{\partial\dot{x}}%
{\partial\sigma}\left(  \sigma_{0}\right)  \in\mu^{1/3}[d_{1},d_{2}%
],\label{eq:sigma0-partials}%
\end{equation}
for constant $d_{1},d_{2}$, $0\notin\left[  d_{1},d_{2}\right]  .$
Stable/unstable manifolds are $S$ symmetric (see (\ref{eq:S-sym}) and Theorem
\ref{lem: transversality-Simo}), hence for a curve coming from the
intersection of the stable manifold with $\{y=0\}$ we shall have same
estimates for $\frac{\partial x}{\partial\sigma},$ and estimates with reversed
sign for $\frac{\partial\dot{x}}{\partial\sigma}$. This by
(\ref{eq:sigma0-partials}) gives a splitting of the two manifolds with angle
of order $\mu^{1/3}.$

\subsection{Derivation of equations for the PRE3BP in rotating coordinates}

Here we derive the Hamiltonian (\ref{eq:H-PRE3BP}) of the PRE3BP in rotating
coordinate system. Let $R(\varphi)$ be the rotation by the angle $\varphi$%
\begin{equation}
R(\varphi)=\left[
\begin{array}
[c]{cc}%
\cos{\varphi} & -\sin{\varphi}\\
\sin{\varphi} & \cos{\varphi}\
\end{array}
\right]  .\nonumber
\end{equation}
The ellipse which is obtained when solving the 2-body problem is given by (see
(\ref{eq:eliptic-orbits}))
\begin{align}
r(t)  &  =\frac{1-e^{2}}{1+e\cos{\psi(t)}}=1-e\cos{\psi(t)}+O(e^{2}%
),\label{eq:r-ellipse}\\
z(t)  &  =r(t)R(\psi(t))\cdot\lbrack1,0]^{T},\nonumber\\
\psi(t)  &  =t+2e\sin{t}+O(e^{2}).\nonumber
\end{align}
The primary with the mass $\mu$ (the planet) has the following location
$z_{p}(t)=(\mu-1)z(t)$, while the primary mass $1-\mu$ (the Sun) is located at
the point $z_{s}(t)=\mu z(t)$.

We shall now compute the distances between the comet and the primaries in
rotating coordinates. Let $r_{1}(t)$ be the square of the distance between the
comet and the Sun, and $r_{2}(t)$ between the comet and the planet. Let
$(x(t),y(t))$ denote the 'rotating' coordinates of the comet and $q(t)$ be the
position of the comet in the 'static' coordinate frame $(x(t),y(t))=R(-t)q(t)$%
. Using the fact the the length of the vector is not changed by the rotation
and (\ref{eq:r-ellipse}) we have
\begin{align*}
r_{1}(t)^{2}=  &  \left\Vert R(t)\cdot(x,y)^{T}-\mu r(t)R(\psi(t))\cdot
\lbrack1,0]^{T}\right\Vert ^{2}\\
=  &  \left\Vert (x,y)^{T}-\mu r(t)R(\psi(t)-t)\cdot\lbrack1,0]^{T}\right\Vert
^{2}\\
=  &  \,x^{2}+y^{2}-2\mu r(t)\left(  x\cos{(\psi(t)-t)}+y\sin{(\psi
(t)-t)}\right)  +\mu^{2}r(t)^{2}\\
=  &  \,x^{2}+y^{2}-2\mu\left(  x\cos{(\psi(t)-t)}+y\sin{(\psi(t)-t)}\right)
+\mu^{2}+\\
&  2\mu e\cos{\psi(t)}\left(  x\cos{(\psi(t)-t)}+y\sin{\psi(t)-t)}\right)
-2\mu^{2}e\cos{\psi(t)}+O(\mu e^{2}).
\end{align*}
Using%
\begin{align*}
\psi(t)  &  =t+2e\sin{t}+O(e^{2}),\\
\cos{(\psi(t))}  &  =\cos{t}-2e\sin^{2}{t}+O(e^{2}),\\
\cos{(\psi(t)-t)}  &  =1+O(e^{2}),\\
\sin{(\psi(t)-t)}  &  =2e\sin{t}+O(e^{2}),
\end{align*}%
\begin{align*}
x\cos{(\psi(t)-t)}+y\sin{(\psi(t)-t)}  &  =x+2ey\sin{t}+O(e^{2}),\\
\cos{(\psi(t))}\left(  x\cos{(\psi(t)-t)}+y\sin{(\psi(t)-t)}\right)   &
=x\cos{t}+O(e),
\end{align*}
in the expression for $r_{1}(t)$ we obtain%
\begin{equation}
r_{1}(t)^{2}=(x-\mu)^{2}+y^{2}+2e\bar{g}(\mu,x,y,t)+O(\mu e^{2}),
\label{eq:D1}%
\end{equation}
where $\bar{g}$ is given in (\ref{eq:g-bar}). Expression for $r_{2}(t)$ is
obtained from (\ref{eq:D1}) with the substitution $\mu\mapsto(\mu-1)$. Observe
that
\begin{equation}
\frac{1}{\sqrt{r^{2}+c}}=\frac{1}{r\sqrt{1+\frac{c}{r^{2}}}}=\frac{1}%
{r}\left(  1-\frac{c}{2r^{2}}+O\left(  \left(  \frac{c}{r^{2}}\right)
^{2}\right)  \right)  . \label{eq:r-expansion-app}%
\end{equation}
We shall use notation $r_{1}$, $r_{2}$ from (\ref{eq:r1-r2-prc3bp}), with
$r_{1}>\delta$ and $r_{2}>\mu^{1/3}\delta$. Note that to apply
(\ref{eq:r-expansion-app}) for $r=r_{2}$ we need to have $e\bar{g}%
(\mu-1,x,y,t)+O(\mu e^{2} )<r_{2}^{2}$, which means that we need to take $e$
sufficiently small so that $e\mu^{-2/3}<\kappa$ for sufficiently small
$\kappa$. Equations (\ref{eq:D1}), (\ref{eq:r-expansion-app}) give
\begin{align}
\frac{1}{r_{1}(t)}  &  =\frac{1}{r_{1}}-\frac{e\bar{g}(\mu,x,y,t)}{r_{1}^{3}%
}+O(\mu^{2}e^{2}),\label{eq:r1-inv}\\
\frac{1}{r_{2}(t)}  &  =\frac{1}{r_{2}}-\frac{e\bar{g}(\mu-1,x,y,t)}{r_{2}%
^{3}}+O(e^{2}\mu^{-5/3}). \label{eq:r2-inv}%
\end{align}
Substituting (\ref{eq:r1-inv}), (\ref{eq:r2-inv}) into
(\ref{eq: H for circular problem}) gives (\ref{eq:H-PRE3BP}).

\end{document}